\documentclass[preprint,12pt,double]{elsarticle}
\usepackage{lineno}
\usepackage{subcaption}
\usepackage{multicol}
\usepackage{float}
\usepackage{siunitx}
\usepackage{nomencl}
\usepackage{amsmath,amssymb,amstext,bm,array,multirow,booktabs,physics} 
\usepackage{mathtools}
\usepackage{graphicx} 
\graphicspath{{figures/}}
\usepackage{appendix}
\usepackage[final]{changes}
\usepackage[capitalize]{cleveref}
\crefname{equation}{Eqn.}{Eqns.}    
\Crefname{equation}{Equation}{Equations}

\DeclarePairedDelimiterX{\inp}[2]{\langle}{\rangle}{#1, #2}

\makenomenclature
\date{25 April, 2019}

\begin{document}


\title{Phase-bounded finite element method for two-fluid incompressible flow systems}

\author[UW]{Tanyakarn Treeratanaphitak}
\author[UW,UWP]{Nasser Mohieddin Abukhdeir\corref{cor}}
\cortext[cor]{Corresponding author}
\ead[url]{http://chemeng.uwaterloo.ca/abukhdeir/}
\ead{nmabukhdeir@uwaterloo.ca}

\address[UW]{Department of Chemical Engineering, University of Waterloo, 200 University Avenue West Waterloo, N2L 3G1, ON, Canada}
\address[UWP]{Department of Physics \& Astronomy, University of Waterloo, 200 University Avenue West Waterloo, N2L 3G1, ON, Canada}

\biboptions{authoryear}
\journal{International Journal of Multiphase Flow}

\begin{frontmatter}

    \begin{abstract}
        An understanding of the hydrodynamics of multiphase processes is essential for their design and operation.
        Multiphase computational fluid dynamics (CFD) simulations enable researchers to gain insight which is inaccessible experimentally.
        The model frequently used to simulate these processes is the two-fluid (Euler-Euler) model where fluids are treated as inter-penetrating continua.
        It is formulated for the multiphase flow regime where one phase is dispersed within another and enables simulation on experimentally relevant scales.
        Phase fractions are used to describe the composition of the mixture and are bounded quantities.
        Consequently, numerical solution methods used in simulations must preserve boundedness for accuracy and physical fidelity.
        In this  work, a numerical method for the two-fluid model is developed in which phase fraction constraints are imposed through the use of an nonlinear variational inequality solver which implicitly imposes inequality constraints.
        The numerical method is verified and compared to an established explicit numerical method.
    \end{abstract}

    \begin{keyword}
        two-phase flow \sep computational fluid dynamics \sep two-fluid model \sep phase fraction boundedness
    \end{keyword}

\end{frontmatter}


\printnomenclature

\section{Introduction}

Two-phase gas-liquid flow systems are prevalent in industrial processes and, consequently, the hydrodynamic behavior of gas-liquid systems is of great interest.
Examples of gas-liquid flow systems include bubble columns \citep{Jakobsen2005,Joshi2001,Rollbusch2015,Shaikh2013,Shah1982}, loop reactors \citep{Becker1994}, cyclones \citep{Erdal1997}, hydrocarbon pipelines \citep{Issa2003,Issa2006} and disengagers \citep{Lane2016a}.
The hydrodynamical behavior of gas-liquid systems is complicated and is dependent on factors such as operating conditions, physical properties, equipment type, etc. \citep{Shah1982, Shaikh2013, Rollbusch2015}.
Conducting experiments to examine the effect of some of these parameters can be difficult as one may be limited to certain experimental measurement techniques and equipment availability \citep{Rollbusch2015}.
Thus the use of experimentation alone to study gas-liquid hydrodynamics is costly, time-consuming and does not easily provide access to the dynamic spatially varying flow field.
Computational fluid dynamics (CFD) simulations are able to augment experimental research in this area through addressing the limitations of current experimental techniques.
Significant progress has been made towards the development of models and numerical methods for the simulation of gas-liquid and other multiphase flows using CFD \citep{Ishii2011,Jakobsen2014,Kolev2011a,Marchisio2007}.
CFD simulations can be used as screening tool for researchers to devise experimental set-ups involving conditions that merit further study.

Current approaches to modeling two-phase flows fall into three main categories: two-fluid (Euler-Euler), Euler-Lagrange and interface-tracking models.
The two-fluid model (TFM) approximates the fluids as inter-penetrating continua with conservation equations formulated for each phase \citep{Ishii2011}.
This method requires interphase momentum transfer constitutive relationships that describe the momentum exchange between each phase.
The selection of these constitutive relationships has been shown to affect the simulation results significantly \citep{Sokolichin2004,Tabib2008}.
This method is the most feasible for industrial-scale simulations as it does not resolve every interface in the system and computing the momentum transfer between every particle/bubble and the fluid, which can be extremely costly in industrial-scale situations.

From a numerical perspective, the solution of the Euler-Euler model presents many challenges including numerical stability \citep{LopezdeBertodano2016,Vaidheeswaran2016a,Vaidheeswaran2016b} and the introduction of inequality constraints \citep{Weller2005,COMSOLCFD}.
The latter is the focus of this work, where constraints on the phase volume fractions, $\alpha_q$, include:
\begin{gather}
    \sum_q \alpha_q = 1,\label{eq:equality_constraint}\\
    0 \le \alpha_q \le 1.\label{eq:inequality_constraint}
\end{gather}
Imposition of \cref{eq:equality_constraint} is straightforward using standard linear or nonlinear solvers.
However, the boundedness constraints on $\{\alpha_i\}$ corresponding to \cref{eq:inequality_constraint} are not and a variety of methods have been developed to impose them including: artificial diffusion \citep{COMSOLCFD}, remapping \citep{Weller2005}, flux limiting \citep{OpenFOAMFoundation2017} and truncation \citep{CFX_Guide}.
All of these approaches either decouple imposition of the inequality constraints from solution of the conservation of mass equations (remapping and truncation) or add \textit{ad hoc} contributions to the conservation of mass equations to impose these inequality constraints.
A significant consequence of decoupled approaches is that local error introduced into the numerical solution is both unquantified and unconstrained. Alternatively, the introduction of \textit{ad hoc} contributions to the conservation of mass equation, such as artificial diffusion, introduce unphysical phenomena into numerical solutions of the model.
Thus the development of a \textit{coupled} numerical method for the solution of TFM would be numerically accurate, with respect to a specified error tolerance, without \textit{ad hoc} modifications of the model.

Within this context, the overall objective of this work is to develop a method for the numerical solution of the TFM which constrains phase fractions through coupled solution of the conservation of mass and imposition of equality (\cref{eq:equality_constraint}) and inequality constraints (\cref{eq:inequality_constraint}).
The numerical solution method proposed involves the use of the finite-element method \citep{Logg2012} and the method of lines \citep{Schiesser1991} in conjunction with a nonlinear solver capable of imposing inequality constraints \citep{petsc-user-ref}.
Specific objectives include:
\begin{enumerate}
  \item Development of a finite-element formulation of the two-fluid model.
  \item Verification of the model with an unbounded version of the model.
  \item Comparison of dispersed flow regime simulations using the numerical method developed in this work and \texttt{OpenFOAM}.
\end{enumerate}

This paper is organized as follows: \cref{sec:background} -- background on the two-fluid model and phase fraction boundedness, \cref{sec:methodology} -- presentation of the model and simulation conditions, \cref{sec:results} -- simulation results on the verification of the model and comparison of the model with \texttt{OpenFOAM} and \cref{sec:conclusions} -- conclusions.

\section{Background} \label{sec:background}

\subsection{Two-Fluid Model} \label{sec:TFM}

In modeling gas-liquid flows using the two-fluid model, each of the phases in the system is considered to be a continuous fluid.
Each of the phases has its own set of conservation equations that are coupled together through interphase transfer terms.
It is impractical to solve for the local instantaneous motion of the fluid, thus averaging schemes are used to solve for the macroscopic flow behavior instead \citep{Ishii2011}.
Time-averaged quantities are denoted by an overbar, $\overline{\centerdot}$.
A double-overbar, $\overline{\overline{\centerdot}}$, denotes a phasic average quantity, which are defined as the time-averaged quantity divided by the phase fraction, $\overline{\centerdot}/\alpha_i$.
A hat, $\widehat{\centerdot}$, denotes a mass-weighted mean phasic average quantity (Favre average), defined as $\overline{\rho_i\centerdot}/\overline{\rho_i}$.

\subsection{Mass Conservation} \label{sec:mass_con}
The general expression for the conservation of mass for a phase $q$, in the absence of interphase mass transfer and reaction, is given as follows \citep{Ishii2011}:\\
\newsavebox\masseq
\savebox\masseq{\vbox{\begin{equation}
        \pdv{(\alpha_{q} \overline{\overline{\rho_{q}}})}{t} + \div (\alpha_{q} \overline{\overline{\rho_{q}}} \widehat{\vb*{v}_{q}})= 0,
        \label{eq:continuity}
        \end{equation}}}\\
\usebox\masseq\\
\nomenclature{$\alpha$}{Time-averaged local phase fraction}
\nomenclature{$\rho$}{Density}
\nomenclature{$t$}{Time}
\nomenclature{$\vb*{v}$}{Velocity vector}
\nomenclature{$\overline{\overline{\centerdot}}$}{Phasic average, $\overline{\centerdot}/\alpha_i$}
\nomenclature{$\overline{\centerdot}$}{Time average (Reynolds average)}
\nomenclature{$\widehat{\centerdot}$}{Mass-weighted mean phasic average (Favre average), $\overline{\rho_i\centerdot}/\overline{\rho_i}$}
where $\alpha_{q}$ is the time-averaged local phase fraction of phase $q$, $\overline{\overline{\rho_{q}}}$ is the time-averaged phasic average density and $\widehat{\vb*{v}_{q}}$ is the time-averaged mass-weighted mean phase velocity.

\subsection{Momentum Conservation} \label{sec:mom_con}

The conservation of momentum for phase $q$ is given as \citep{Ishii2011}:
\begin{equation}
\begin{split}
\pdv{\qty(\alpha_{q} \overline{\overline{\rho_{q}}} \widehat{\vb*{v}_{q}})}{t} + \div (\alpha_{q} \overline{\overline{\rho_{q}}} \widehat{\vb*{v}_{q}} \widehat{\vb*{v}_{q}}) &= -\grad\qty(\alpha_q \overline{\overline{P_{q}}}) + \div \qty(\alpha_{q} {\overline{\overline{\vb*{\tau}_{q}}}})+\alpha_{q} \overline{\overline{\rho_{q}}} {\widehat{\vb*{g}_{q}}} + {\vb*{M}_{q}} \\
&\qquad + \overline{\overline{P_{q,i}}} \grad \alpha_q - \grad\alpha_q\vdot \overline{\overline{\vb*{\tau}_{q,i}}},
\end{split}
\label{eq:momentum_eq}
\end{equation}
\nomenclature{$P$}{Pressure}
\nomenclature{$\vb*{\tau}$}{Stress tensor}
\nomenclature{$\vb*{g}$}{Gravitational acceleration}
\nomenclature{$\vb*{M}$}{Momentum source term}
where $\overline{\overline{P_{q}}}$ is the time-averaged phasic pressure, ${\overline{\overline{\vb*{\tau}_{q}}}}$ is the time-averaged phasic viscous stress tensor, ${\widehat{\vb*{g}_{q}}}$ is the time-averaged mass-weighted mean phase gravitational acceleration, ${\vb*{M}_{q}}$ is the interphase momentum source term, $\overline{\overline{P_{q,i}}} \grad \alpha_q$ and $\grad\alpha_q \vdot \overline{\overline{\vb*{\tau}_{q,i}}}$ are the contributions of interfacial stresses and the subscripts indicate the fluid-dispersed phase pairing and type of force, respectively.
In the dispersed flow regime, the interfacial pressure and shear stress of the continuous, $c$, and dispersed, $d$, phases can be assumed to be equal to each other $\overline{\overline{P_{c,i}}} \approx \overline{\overline{P_{d,i}}} = \overline{\overline{P_{int}}}$ and $\overline{\overline{\vb*{\tau}_{c,i}}} \approx \overline{\overline{\vb*{\tau}_{d,i}}}$ \citep{Drew1998,Ishii2011}.
Additionally, the pressure of the dispersed phase can be approximated by the interfacial pressure, $\overline{\overline{P_{d}}} \approx \overline{\overline{P_{int}}}$ \citep{Ishii2011}:
\begin{subequations}
    \begin{align}
    \begin{split}
    \pdv{\qty(\alpha_{c} \overline{\overline{\rho_{c}}} \widehat{\vb*{v}_{c}})}{t} + \div (\alpha_{c} \overline{\overline{\rho_{c}}} \widehat{\vb*{v}_{c}} \widehat{\vb*{v}_{c}}) &= -\alpha_c\grad \overline{\overline{P_{c}}} + \div \qty(\alpha_{c} {\overline{\overline{\vb*{\tau}_{c}}}})+\alpha_{c} \overline{\overline{\rho_{c}}} {\widehat{\vb*{g}}} + {\vb*{M}_{c}} \\
    &\qquad + \qty(\overline{\overline{P_{int}}} - \overline{\overline{P_{c}}})\grad \alpha_c - \grad\alpha_q\vdot \overline{\overline{\vb*{\tau}_{c,i}}},
    \end{split}\label{eq:momentum_eq_c} \\
    \begin{split}
    \pdv{\qty(\alpha_{d} \overline{\overline{\rho_{d}}} \widehat{\vb*{v}_{d}})}{t} + \div (\alpha_{d} \overline{\overline{\rho_{d}}} \widehat{\vb*{v}_{d}} \widehat{\vb*{v}_{d}}) &= -\alpha_d\grad \overline{\overline{P_{int}}} + \div \qty(\alpha_{d} {\overline{\overline{\vb*{\tau}_{d}}}})+\alpha_{d} \overline{\overline{\rho_{d}}} {\widehat{\vb*{g}}} + {\vb*{M}_{d}} \\
    &\qquad - \grad\alpha_d\vdot \overline{\overline{\vb*{\tau}_{c,i}}}.
    \end{split}
    \label{eq:momentum_eq_d}
    \end{align}%
\end{subequations}
The effect of the interfacial shear stress is significant in the segregated flow regime \citep{Ishii1984}.
In this work, the focus will be on the dispersed flow regime, therefore, the interfacial shear stress contribution assumed to be negligible.
The conservation of momentum for a dispersed flow system is thus:
\newsavebox\momentumeq
\savebox\momentumeq{\vbox{\begin{subequations}
            \begin{align}
            \begin{split}
            \pdv{\qty(\alpha_{c} \overline{\overline{\rho_{c}}} \widehat{\vb*{v}_{c}})}{t} + \div (\alpha_{c} \overline{\overline{\rho_{c}}} \widehat{\vb*{v}_{c}} \widehat{\vb*{v}_{c}}) &= -\alpha_c\grad \overline{\overline{P_{c}}} + \div \qty(\alpha_{c} {\overline{\overline{\vb*{\tau}_{c}}}})+\alpha_{c} \overline{\overline{\rho_{c}}} {\widehat{\vb*{g}}} + {\vb*{M}_{c}} \\
            &\qquad + \qty(\overline{\overline{P_{int}}} - \overline{\overline{P_{c}}})\grad \alpha_c,
            \end{split}\label{eq:momentum_eq_c2} \\
            \pdv{\qty(\alpha_{d} \overline{\overline{\rho_{d}}} \widehat{\vb*{v}_{d}})}{t} + \div (\alpha_{d} \overline{\overline{\rho_{d}}} \widehat{\vb*{v}_{d}} \widehat{\vb*{v}_{d}}) &= -\alpha_d\grad \overline{\overline{P_{int}}} + \div \qty(\alpha_{d} {\overline{\overline{\vb*{\tau}_{d}}}})+\alpha_{d} \overline{\overline{\rho_{d}}} {\widehat{\vb*{g}}} + {\vb*{M}_{d}}.
            \label{eq:momentum_eq_d2}
            \end{align}\label{eq:momentum_eq_2}%
\end{subequations}}}\\
\usebox\momentumeq\\
For the sake of brevity, the average notations are dropped in the subsequent sections.
The quantities are assumed to be averaged quantities corresponding to the definition above.

\subsection{Interfacial Forces} \label{sec:interfacial_forces}

The interphase momentum transfer term ${\vb*{M}_{q}}$, defined as the transfer of momentum \emph{into} phase $q$, is the sum of the contributions from the different modes of momentum transfer: drag, lift, virtual mass, wall lubrication, \textit{etc}.
The momentum transfer term for the continuous phase with the aforementioned contributions is \citep{Ishii2011,Lahey2001,Antal1991}:
\begin{equation}
{\vb*{M}_{c}} = {\vb*{M}_{c,drag}}+{\vb*{M}_{c,lift}}  + {\vb*{M}_{c,virtual \ mass}} + \vb*{M}_{c, wall}.
\label{eq:mom_trans}
\end{equation}
Given that the momentum exchange between the phases should sum to zero, momentum transfer of the dispersed phase $d$ is given as:
\begin{equation}
{\vb*{M}_{c}} = -{\vb*{M}_{d}}.
\label{eq:mom_trans3}
\end{equation}

\subsubsection{Drag Force}\label{sec:drag_force}

The drag force term in \cref{eq:mom_trans} is the sum of the form and skin drag forces which are due to the imbalance of pressure and shear forces at the interface, respectively \citep{Ishii2011}.
The interphase momentum transfer of phase $c$ due to drag for a dispersed phase $d$ in fluid $c$, ${\vb*{M}_{c,drag}}$, is given as \citep{Ishii2011}:
\begin{equation}
{\vb*{M}_{c,drag}} = \frac{1}{2}\rho_c\alpha_d \frac{C_D}{r_d}\norm{\vb*{v}_r}\vb*{v}_r,
\label{eq:drag_force}
\end{equation}
\nomenclature{$C_D$}{Drag coefficient}%
\nomenclature{$d$}{Diameter}%
\nomenclature{$r_d$}{Volume to projected area ratio}%
where $r_d$ is the ratio of the volume to projected area of the bubble/particle, $C_D$ is the drag coefficient and $\vb*{v}_r$ is the relative velocity between the dispersed and continuous phases, $\vb*{v}_r = \vb*{v}_d - \vb*{v}_c$.
For spherical bubbles/particles,
\begin{equation}
{\vb*{M}_{c,drag}} =  \frac{3}{4}\rho_c\alpha_d \frac{C_D}{d_d}\norm{\vb*{v}_r}\vb*{v}_r,
\label{eq:drag_spherical}
\end{equation}
where $d_d$ is the bubble/particle diameter.

\subsubsection{Lift Force}

Lift is the force that governs the transverse movement of the dispersed phase in a fluid and is a result of shear forces and the asymmetric pressure distribution around the dispersed particle/bubble \citep{Sokolichin2004,Tomiyama2002, Drew1987}.
The direction of lift is perpendicular to the direction of flow.
The expression for the momentum transfer to the continuous fluid $c$ due to lift is given as \citep{Drew1987,Drew1993}:
\begin{equation}
{\vb*{M}_{c,lift}} = C_L \rho_{c} \alpha_{d} \vb*{v}_r \times \qty(\curl \vb*{v}_c),
\label{eq:lift_drew-lahey}
\end{equation}
\nomenclature{$C_L$}{Lift coefficient}
where $C_L$ is the lift coefficient.

\subsubsection{Virtual Mass Force}

The virtual mass (or added mass) force is related to the acceleration of one phase in relation to another \citep{Drew1987}.
The momentum transfer to phase $c$ due to virtual mass force is given as \citep{Drew1987,Drew1993}:
\begin{equation}
{\vb*{M}_{c,virtual \ mass}} = \alpha_d \rho_c C_{VM}\qty(\frac{D_d \vb*{v}_d}{Dt} - \frac{D_c \vb*{v}_c}{Dt} ),\label{eq:virtual_mass}
\end{equation}
\nomenclature{$\frac{D_{i} \centerdot}{Dt}$}{Material derivative, $\pdv{\centerdot}{t} + \vb*{v}_i \vdot\grad \centerdot$}
\nomenclature{$C_{VM}$}{Virtual mass coefficient}
where $C_{VM}$ is the virtual mass coefficient, $D_d/Dt$ and $D_c/Dt$ are the material derivatives with respect to phases $d$ and $c$, respectively.

\subsubsection{Wall Lubrication Force}

The wall lubrication force is a wall effect that occurs in bubbly flow where the continuous phase wets the walls.
It occurs when a bubble's proximity to the wall results in asymmetric drainage of fluid around the bubble.
The side that is close to the wall will drain slower due to the no-slip condition.
The asymmetry creates a hydrodynamic force normal to the wall that pushes the bubble away from the wall \citep{Antal1991}.
The wall lubrication force is given as \citep{Frank2008}:
\begin{equation}
\vb*{M}_{c, wall} = -C_W \alpha_d \rho_c \norm{\vb*{v}_r - \qty(\vb*{v}_r\vdot \vb*{n}_W)\vb*{n}_W}^2 \vb*{n}_W,
\end{equation}
\nomenclature{$C_W$}{Wall coefficient}
where $C_W$ is the wall coefficient and $\vb*{n}_W$ is the unit normal outward on the wall.
As the distance between the bubble and the wall increases, $C_W$ will taper off towards zero.

\subsubsection{Interfacial Pressure}

The interfacial pressure is determined from a volume average of the solution of potential flow around a single sphere \citep{Stuhmiller1977,Antal1991}.
This interfacial pressure is given by:
\begin{equation}
P_{c,i} = P_c - C_P \rho_c \vb*{v}_r \vdot \vb*{v}_r,
\label{eq:interfacial_pressure}
\end{equation}
\nomenclature{$C_{P}$}{Interfacial pressure coefficient}
where $C_P$ is the interfacial pressure coefficient.
For the case where the particle/bubble size distribution is uniform, $C_P = 0.25$ \citep{Pauchon1986}.

\subsection{Phase Fraction Boundedness}

An important aspect of the two-fluid model is the boundedness of phase fractions.
\Cref{eq:equality_constraint} is inherently satisfied in the two-fluid model when solving the conservation of mass equation for only one of the phases, typically the dispersed phase.
The phase fraction of the continuous phase is simply $1-\alpha_d$.
\Cref{eq:inequality_constraint} is not inherently satisfied and additional steps must be taken to ensure that the individual phase fractions are bounded.

The two-fluid model has been implemented in several commercial software packages \citep{CFX_Guide,Fluent_Guide,COMSOLCFD} that are readily available.
Other implementations such as \texttt{NEPTUNE\_CFD} \citep{Bartosiewicz2008,Mimouni2008} and \texttt{OpenFOAM} \citep{Weller1998,Weller2005} also exist.
A comparative study of some of the TFM implementations is available in the literature \citep{Rzehak2015,Bartosiewicz2008}.
Some of the packages \citep{Fluent_Guide,COMSOLCFD,Weller2005,Damian2013} include details on how phase fraction boundedness is enforced.

The different approaches to maintaining phase fraction boundedness available in literature can be broken down into four categories: thresholding \citep{Fluent_User_Guide}, flux limiting \citep{OpenFOAM}, artificial diffusion \citep{COMSOLCFD} and remapping \citep{Oliveira2003,Weller2005}.
The use of phase fraction bounding techniques can have an effect on the fidelity of the solution.
In this section, existing approaches in maintaining phase fraction boundedness are reviewed.

Depending on the formulation of the momentum equation used, thresholding can be necessary to avoid issues with division by zero.
The value of the phase fraction below a threshold is set to be equal to the threshold, typically a small value.
One can also choose to threshold the phase fraction outside the momentum equation, this involves zeroing out negative phase fraction values and setting all values above one to be equal to one after the solving \cref{eq:continuity}.
This is the approach taken in \citet{CFX_Guide}.
Thresholding the phase fraction from \cref{eq:continuity} will result in a change in the profile and gradient of the phase fraction.

\citet{Oliveira2003} implemented a two-equation method to satisfy \cref{eq:equality_constraint,eq:inequality_constraint} in their in-house code.
\Cref{eq:continuity} was discretized using an upwind scheme and solved for both the dispersed and continuous phases.
The upwinding ensured that the lower bound in \cref{eq:inequality_constraint} is satisfied.
The resulting phase fractions were then rescaled by a factor of $1/\qty(\alpha_d^*+\alpha_c^*)$, where $\alpha^*_d$ and $\alpha_c^*$ are the values obtained from the conservation of mass equation, satisfying \cref{eq:equality_constraint}.
Since the rescaled individual phase fractions satisfy the lower inequality bound and \cref{eq:equality_constraint}, the upper inequality bound is also satisfied.

In previous versions of \texttt{OpenFOAM}, \citet{Weller2005} reformulated the conservation of mass equation such that the phase fraction can be bounded when the conservative form is used.
The dispersed phase velocity $\vb*{v}_d$ was decomposed into mean and relative components \citep{Weller2005}:
\begin{equation}
\vb*{v}_d = \alpha_c \vb*{v}_c + \alpha_d \vb*{v}_d + \alpha_c \vb*{v}_r. \label{eq:weller_vd}
\end{equation}
Substituting \cref{eq:weller_vd} into \cref{eq:continuity} and dividing by $\rho_d$:
\begin{equation}
\pdv{\alpha_d}{t} + \div{\qty[\alpha_d \qty(\alpha_c \vb*{v}_c + \alpha_d \vb*{v}_d)]} + \div{\qty(\vb*{v}_r \alpha_d \alpha_c)} = 0. \label{eq:weller_continuity}
\end{equation}
The resulting equation is nonlinear and the boundedness of the solution may be compromised when using a higher order spatial discretization scheme \citep{Rusche2002}.
To solve the equation using an iterative linear solver while maintaining the boundedness of the phase fractions, the phase fraction was remapped using a quadratic equation that is a function of both phase fractions \citep{Weller2005}:
\begin{equation}
\alpha_d^* = \frac{1}{2}\qty[1-\qty(1-\alpha_d)^2 + \qty(1-\alpha_c)^2].
\end{equation}

In more recent versions of \texttt{OpenFOAM}, phase fraction boundedness is ensured through the use of a limiter that is based on flux corrected transport called multidimensional universal limiter for explicit solution (MULES) \citep{OpenFOAM}.
MULES allows for the possibility of specifying global minimum and maximum values for a given field \citep{Damian2013}. 
\added{It uses the flux-corrected transport (FCT) approach, which is specific to purely conservative PDEs with inequality constraints.
The method reduces the mass-flux into or out of a cell when that flux would violate an inequality constraint.
The conservation of mass equation in the TFM can be represented in divergence form, where the mass flux between cells results solely from convection.
Thus, MULES constrains the velocity component parallel to the face of cells in which the current iteration of the momentum solve would result in an out of bounds phase fraction.
This method is compatible with explicit solution methods and convenient for use with the finite volume method in that fluxes between cells are already computed.
The drawbacks of this approach include the requirement that fluxes between cells must be computed, closed-form expressions relating these fluxes to constrained variables must be derived, and in \texttt{OpenFOAM} the FCT method is applied during the explicit corrector iterations of the momentum solve.
The consequence of this last drawback is that modification of the corrected velocity field resulting from the MULES method could alter its divergence-free property.}

The artificial diffusion approach adds a diffusion term to the conservation of mass equation that will help regularize the solution.
In the commercial package \texttt{COMSOL} \citep{COMSOLCFD}, artificial diffusion is used in the conservation of mass equation of the dispersed phase to minimize the possibility of a negative phase fraction:
\begin{equation}
\pdv{\alpha_d}{t} + \div{\qty(\alpha_d \vb*{v}_d)} = -\div{\qty(-\nu_b \grad{\alpha_d})},
\end{equation}
\nomenclature{$\mu$}{Dynamic viscosity}%
\nomenclature{$\nu$}{Kinematic viscosity}%
where $\nu_b$ is the ``barrier'' viscosity:
\begin{equation}
\nu_b = \frac{\mu_d}{\rho_d}\qty(\exp\qty[\max\qty(-\frac{\alpha_d}{0.0025},0)]-1).
\end{equation}
The barrier viscosity is nonzero when the phase fraction is a negative value.
The artificial diffusion term only minimizes the possibility of a negative phase fraction, it does not guarantee that the phase fraction will be positive.
Artificial diffusion can alter the governing equation even if the diffusion term is only active when the phase fraction is negative.
When the conservation of momentum equation is solved, the phase fractions are thresholded to be between zero and one to regularize the solution \citep{COMSOLCFD}.
In \texttt{COMSOL}, the convective form of the momentum equation is scaled by the phase fraction, which results in $1/\alpha$ terms in the momentum equation, these terms are also thresholded.

\section{Methodology}\label{sec:methodology}

The conservation of momentum equation defined in \cref{eq:momentum_eq_2} is conservative, but its solution becomes degenerate as the phase fraction approaches zero.
To avoid this issue, the dimensionless convective form of the momentum equation, scaled by $\alpha_q\rho_q$, is used instead.
The contribution of the interfacial shear stress is neglected.
In this study, the fluid-fluid system of interest is a gas-liquid system where liquid, $l$, is the continuous phase and gas, $g$, is the dispersed phase, the governing equations will reflect this from hereon in.
The dimensionless parameters are defined as follow: $\tilde{\boldsymbol{v}} = \boldsymbol{v}/v_s$, $\tilde{t} = t/t_s$, $\tilde{\boldsymbol{x}} = \boldsymbol{x}/x_s$, $\tilde{P} = (P-P_0)/P_s$,  $\tilde{\boldsymbol{g}} = \boldsymbol{g}/g_s$, $\tilde{\nabla} = x_s\nabla$ and $ \tilde{d_b} = d_b/x_s$.

\begin{table}
    \caption{Dimensionless groups}
    \label{tab:parameters}
    \centering
    \begin{tabular}{lc}
        \hline
        Parameter                   &          Expression           \\ \hline
        Time                        &        $t_s = v_s/x_s$        \\
        \multirow{2}{5em}{Pressure} &     $P_s = \rho_l g_s h$      \\
                                    &           $P_0 = 0$           \\
        Euler number                &   $Eu_q = P_s/\rho_q v_s^2$   \\
        Reynolds number             & $Re_q = \rho_q v_s x_s/\mu_q$ \\
        Froude number               &   $Fr = v_s/\sqrt{g_s x_s}$   \\ \hline
    \end{tabular}
\end{table}
\nomenclature{$Re$}{Reynolds number}%
\nomenclature{$Eu$}{Euler number}%
\nomenclature{$Fr$}{Froude number}%
This results in the following scaled equations:
\begin{subequations}
    \begin{align}
    \begin{split}  \pdv{\tilde{\vb*{v}}_{l}}{ \tilde{t}} +  \tilde{\vb*{v}}_{l}\vdot\tilde{\grad} \tilde{\vb*{v}}_{l} &= -Eu_l \tilde{\grad} \tilde{P}_l + \frac{1}{Re_l}\frac{\tilde{\grad}\alpha_l \vdot \tilde{\vb*{\tau}}_l}{\alpha_l}  + \frac{1}{Re_l}\tilde{\div}\tilde{\vb*{\tau}}_l+ \frac{1}{Fr^2} \tilde{\vb*{g}} \\
    & \quad +\frac{3}{4}\frac{\alpha_g}{\alpha_l} \frac{C_{D}}{\tilde{d}_b}\norm{\tilde{\vb*{v}}_r}\tilde{\vb*{v}}_r - C_P \tilde{\vb*{v}}_r \vdot\tilde{ \vb*{v}}_r\frac{\tilde{\grad}\alpha_l}{\alpha_l} , \end{split} \label{scaled_mom_l}\\
    \begin{split} \pdv{\tilde{\vb*{v}}_{g}}{\tilde{t}} +   \tilde{\vb*{v}}_{g}\vdot\tilde{\grad} \tilde{\vb*{v}}_{g} &= -Eu_g \tilde{\grad} \qty(\tilde{P}_l - C_P \tilde{\vb*{v}}_r \vdot\tilde{ \vb*{v}}_r\frac{\rho_l}{\rho_g}) + \frac{1}{Re_g}\frac{\tilde{\grad}\alpha_g \vdot \tilde{\vb*{\tau}}_g}{\alpha_g} \\
    & \quad + \frac{1}{Re_g}\tilde{\div}\tilde{\vb*{\tau}}_g+ \frac{1}{Fr^2}  \tilde{\vb*{g}}  -\frac{3}{4}\frac{\rho_l}{\rho_g}\frac{C_{D}}{\tilde{d}_b}\norm{\tilde{\vb*{v}}_r}\tilde{\vb*{v}}_r , \end{split} \label{eq:scaled_mom_g}\\
    \pdv{ \alpha_g}{\tilde{t}} + \tilde{\div}\qty(\alpha_g\tilde{\vb*{v}}_g) &= 0, \label{eq:scaled_mass}\\
    \alpha_l &= 1 - \alpha_g,\label{eq:alpha_l}
    \end{align}%
\end{subequations}
where the dimensionless groups are given in \cref{tab:parameters} and \cref{eq:alpha_l} is a result of \cref{eq:equality_constraint}.

\subsection{Modified Incremental Pressure Correction Scheme}\label{sec:IPCS-decoupled}

The governing equations are solved using a modification of the incremental pressure correction scheme (IPCS) \citep{Goda1979}.
IPCS was originally developed for single-phase flow and has been shown to be efficient and accurate \citep{Logg2012}. 
In this work, a modification is made to accommodate for the addition of the conservation of mass and second conservation of momentum equation in the two-fluid model.
Explicit Euler time discretization is used to demonstrate the form of the discretized equations.
There are four steps in the scheme, the first step is to compute a tentative velocity, $\tilde{\vb*{v}}_q^*$, using pressure and velocities from the previous time step:
\begin{subequations}
    \begin{align}
    \begin{split}
    \frac{\tilde{\vb*{v}}_{l}^{*}- \tilde{\vb*{v}}_{l}^n}{\Delta t} +  \tilde{\vb*{v}}_{l}^n\vdot\tilde{\grad} \tilde{\vb*{v}}_{l}^n &= -Eu_l \tilde{\grad} \tilde{P}^n_l + \frac{1}{Re_l}\frac{\tilde{\grad}\alpha_l^n \vdot \tilde{\vb*{\tau}}_l^{n+\frac{1}{2}}}{\alpha_l^n} + \frac{1}{Re_l}\tilde{\div}\tilde{\vb*{\tau}}_l^{n+\frac{1}{2}}\\
    & \quad+ \frac{1}{Fr^2} \tilde{\vb*{g}}  + \frac{3}{4}\frac{\alpha_g^n}{\alpha_l^n} \frac{C_{D}}{\tilde{d}_b}\norm{\tilde{\vb*{v}}_r^n}\tilde{\vb*{v}}_r^n  - C_P \tilde{\vb*{v}}_r^n \vdot\tilde{ \vb*{v}}_r^n\frac{\tilde{\grad}\alpha_l^n}{\alpha_l^n}
    \end{split}& & \text{in} \quad\Omega,\label{eq:mom_l_discretized}\\
    \begin{split}
    \frac{\tilde{\vb*{v}}_{g}^{*}- \tilde{\vb*{v}}_{g}^n}{\Delta t} + \tilde{\vb*{v}}_{g}^n\vdot\tilde{\grad} \tilde{\vb*{v}}_{g}^n &= -Eu_g \tilde{\grad} \qty(\tilde{P}_l^n - C_P \tilde{\vb*{v}}_r^n \vdot\tilde{ \vb*{v}}_r^n\frac{\rho_l}{\rho_g})\\
    & \quad + \frac{1}{Re_g}\frac{\tilde{\grad}\alpha_g^n \vdot \tilde{\vb*{\tau}}_g^{n+\frac{1}{2}}}{\alpha_g^n} + \frac{1}{Re_g}\tilde{\div}\tilde{\vb*{\tau}}_g^{n+\frac{1}{2}} + \frac{1}{Fr^2}  \tilde{\vb*{g}}\\
    & \quad -\frac{3}{4}\frac{\rho_l}{\rho_g}\frac{C_{D}}{\tilde{d}_b}\norm{\tilde{\vb*{v}}_r^n}\tilde{\vb*{v}}_r^n
    \end{split}& & \text{in} \quad\Omega,\label{eq:mom_g_discretized}\\
    \tilde{\vb*{v}}_q^* &= \tilde{\vb*{v}}_{q,BC}^{n+1} & & \text{on} \quad\Gamma_D, \\
    \vb*{n}\vdot\tilde{\vb*{\tau}}_q^{n+\frac{1}{2}} &= \vb*{0} & & \text{on} \quad\Gamma_N, \label{eq:mom_NBC}
    \end{align}\label{eq:TFM_momentum}%
\end{subequations}
\nomenclature{$\Omega$}{Simulation domain}%
\nomenclature{$\Gamma$}{Simulation domain boundary}%
\nomenclature{$\Delta t$}{Time step}%
where $\tilde{\vb*{\tau}}_q^{n+\frac{1}{2}} = \qty(\tilde{\vb*{\tau}}_q^{n+1}+\tilde{\vb*{\tau}}_q^{n})/2$.
The $\qty(\tilde{\grad}{\alpha^n_q})/\alpha^n_q$ term is approximated by $\tilde{\grad}\qty(\ln \alpha'_q)$ where $\alpha'_q$ is $\alpha^n_q$ thresholded to be above a minimum value, $10^{-5}$ in this work.
Then, the tentative velocity is used to compute an update to the pressure field.
Using the conservative form of the momentum equations and the incompressibility criterion of the two-fluid model, $\div\sum_q \alpha_q\vb*{v}_q = 0$, the following pressure Poisson equation is obtained:
\begin{subequations}
    \begin{align}
    \div \qty[\sum_q Eu_q \alpha_q^n \grad (\tilde{P}_l^{n+1}-\tilde{P}_l^n)] &=  \div \sum_q \qty(\frac{\alpha_q^n\tilde{\vb*{v}}_q^{*}}{\Delta t}) & & \text{in} \quad\Omega, \label{eq:pp_discretized}\\
    \vb*{n}\vdot\grad \qty(\tilde{P}_l^{n+1}-\tilde{P}_l^n) &= 0 & & \text{on} \quad\Gamma_D, \label{eq:pp_NBC} \\
    \tilde{P}_l^{n+1} &= \tilde{P}^{n+1}_{l,BC} & & \text{on} \quad\Gamma_N.
    \end{align}%
\end{subequations}
This pressure is used to update the velocity fields:
\begin{subequations}
    \begin{align}
    \frac{\tilde{\vb*{v}}_{l}^{n+1}- \tilde{\vb*{v}}_{l}^*}{\Delta t} &=  -Eu_l \grad \qty(\tilde{P}^{n+1}_l-\tilde{P}^n_l) & & \text{in} \quad\Omega, \label{eq:vl_update_discretized}\\
    \frac{\tilde{\vb*{v}}_{g}^{n+1}- \tilde{\vb*{v}}_{g}^*}{\Delta t} &=  -Eu_g \grad \qty(\tilde{P}^{n+1}_l-\tilde{P}^n_l) & & \text{in} \quad\Omega, \label{eq:vg_update_discretized}\\
    \tilde{\vb*{v}}_q^{n+1} &= \tilde{\vb*{v}}_{q,BC}^{n+1} & & \text{on} \quad\Gamma_D.
    \end{align}%
\end{subequations}
Finally, the phase fractions are computed using the updated velocity field:
\begin{subequations}
    \begin{align}
    \frac{\alpha_g^{n+1} - \alpha_g^n}{\Delta t} + \tilde{\div}\qty(\alpha_g^{n+1}\tilde{\vb*{v}}_g^{n+1}) &= 0 & & \text{in} \quad\Omega,\label{eq:alpha_update_discretized}\\
    \alpha_g^{n+1} &= \alpha_{g,BC}^{n+1} & & \text{on} \quad\Gamma_D,\\
    \alpha_l^{n+1} &= 1 - \alpha_g^{n+1} & & \text{in} \quad\Omega.
    \end{align}%
\end{subequations}

\subsection{Simulation Conditions}\label{sec:conditions}

The two-fluid model is generally well-suited for dilute systems where the phase fraction of the dispersed phase is less than 3\% \citep{Behzadi2004}.
At higher phase fractions, factors such as turbulence, swarming, \emph{etc.} play an non-negligible role in the hydrodynamical behavior of the system, resulting in interphase momentum transfer terms that are dependent on the flow regime.
In this study, the geometry and physical properties of the system are chosen such that the flow remains dispersed for a long period of time to avoid such dependencies.

\begin{figure}
    \centering
    \includegraphics[width=0.4\textwidth]{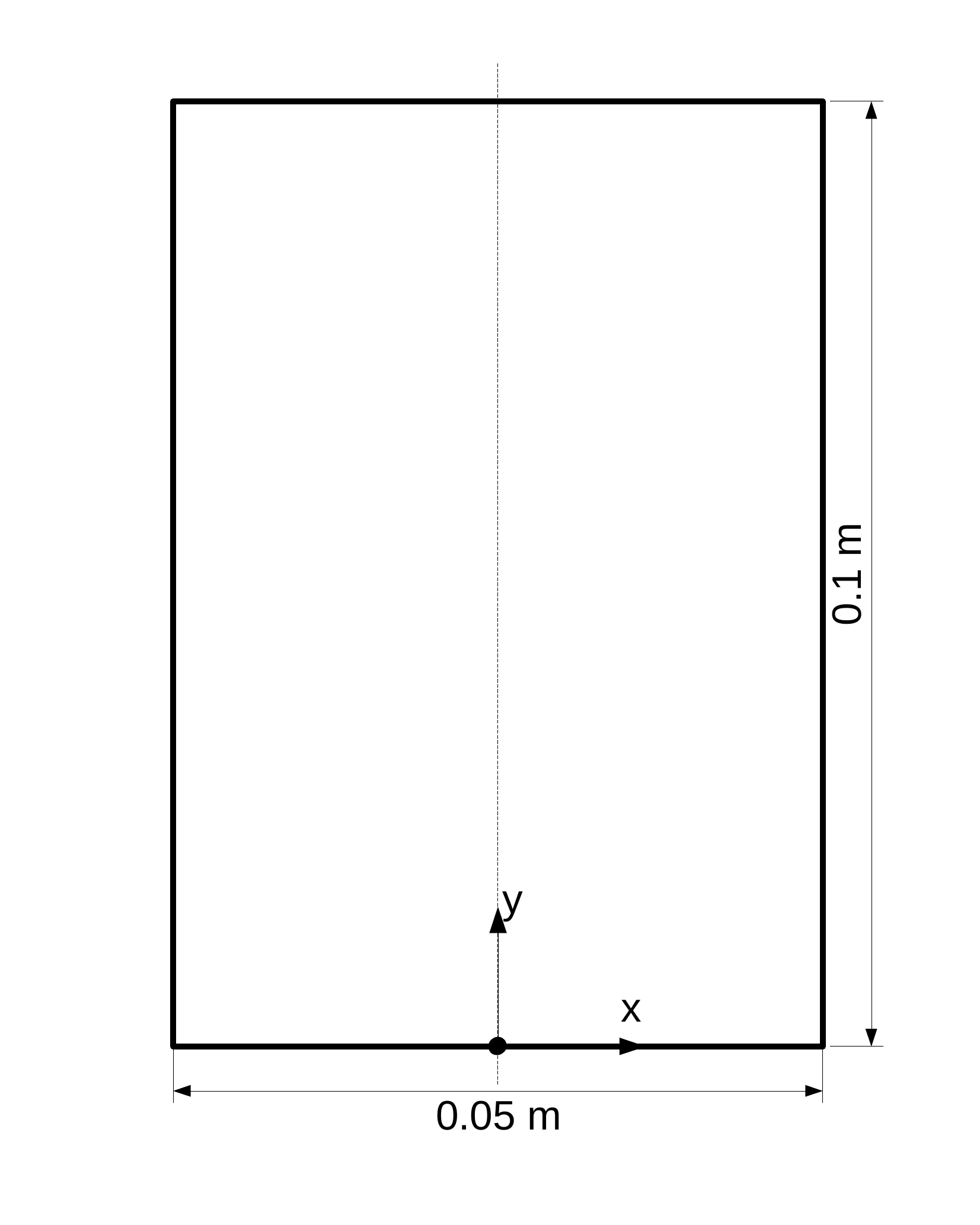}
    \caption{Simulation domain}
    \label{fig:domain}
\end{figure}

The simulation domain used is a two-dimensional channel with gas phase injected from the bottom, as shown in \cref{fig:domain}.
The interface drag coefficient $C_D$ is approximated using the Schiller-Naumann drag expression \citep{Schiller1935}.
Momentum transfer due to lift and virtual mass are neglected.
The Schiller-Naumann drag expression, physical properties and initial and boundary conditions are summarized in \cref{tab:phys_prop,,tab:bcs_var_CD}.
The inlet gas velocity and phase fraction profiles follow a Gaussian distribution and free-slip boundary condition for the gas phase at the walls is used.
In lieu of a wall lubrication force, the gas phase fraction at the walls is set to zero (liquid wets the wall).

\begin{table}
    \centering
    \caption{Physical properties}
    \label{tab:phys_prop}
    \begin{tabular}{lc}
        \hline
        Property & \\
        \hline
        Gas density ($\si{\kg\per\meter^3}$) & $10$ \\
        Liquid density ($\si{\kg\per\meter^3}$) & $1000$ \\
        Gas viscosity ($\si{\pascal \second}$) & $2\times 10^{-5}$ \\
        Liquid viscosity ($\si{\pascal \second}$) & $5\times 10^{-3}$ \\
        Bubble diameter ($\si{\meter}$) & $10^{-3}$ \\
        Drag constant & $\max\qty[\frac{24}{Re}\qty(1+0.15Re^{0.687}),0.44], Re = \frac{\rho_l \norm{\vb*{v}_r}d_b}{\mu_l}$ \\
        \hline
    \end{tabular}
\end{table}

\begin{table}
    \centering
    \caption{Initial and boundary conditions.}
    \label{tab:bcs_var_CD}
    \begin{tabular}{lc}
        \hline
        & Condition \\
        \hline
        \multirow{3}{*}{Initial} & $\alpha_g\qty(\vb*{x},0) = 0$ \\
        & $\vb*{v}_g\qty(\vb*{x},0) = \vb*{v}_l\qty(\vb*{x},0) = \vb{0}$ \\
        & $P\qty(\vb*{x}, 0) = \rho_l g_s (0.1 - y)$ \\
        \hline
        \multirow{4}{*}{Inlet} & $\vb*{v}_g\qty(x,0,t) = \qty(0, \min\qty(\frac{t}{t_0}, 1)0.0616\exp\qty[-\frac{\qty(\frac{x}{0.025})^2}{2\sigma^2}]), t_0 = 0.625~\si{\second}, \sigma = 0.1$ \\
        & $\vb*{v}_l\qty(x,0,0) = \vb*{0} $ \\
        &  $\alpha_{g}\qty(x,0,t) = \min\qty(\frac{t}{t_0}, 1)0.026\exp\qty[-\frac{\qty(\frac{x}{0.025})^2}{2\sigma^2}], t_0 = 0.625~\si{\second}, \sigma = 0.1$ \\
        & $\vb*{n}\vdot\grad{\qty(P_l\qty(x,0,t) - P_l\qty(x,0,t-\Delta t))} = 0$ \\
        \hline
        \multirow{3}{*}{Walls} & $\vb*{v}_l\qty(\pm 0.025,y,t) = \vb{0}$ \\
        & $\vb*{n}\vdot\vb*{v}_g\qty(\pm 0.025,y,t) = 0$ \\
        & $\vb*{n}\vdot\grad{\vb*{v}_g} \vdot \vb*{t}\qty(\pm 0.025,y,t) = 0 $\\
        & $\alpha_g \qty(\pm 0.025,y,t) = 0 $ \\
        & $\vb*{n}\vdot\grad{\qty(P_l\qty(\pm 0.025,y,t) - P_l\qty(\pm 0.025,y,t - \Delta t))} = 0$ \\
        \hline
        \multirow{3}{*}{Outlet} & $\vb*{n}\vdot\vb*{\tau}_g\qty(x,0.1,t)\vdot\vb*{n} = \vb*{n}\vdot\vb*{\tau}_l\qty(x,0.1,t)\vdot\vb*{n} = 0$ \\
        & $\vb*{t}\vdot\vb*{v}_g\qty(x,0.1,t) = \vb*{t}\vdot\vb*{v}_l\qty(x,0.1,t) = 0$ \\
        & $P_l\qty(x,0.1,t) = 0$ \\
        \hline
    \end{tabular}
\end{table}

\subsection{Numerical Methods}

Simulations with the conditions described in the previous section are carried out using the IPCS solver also presented in the previous section.
The equations are solved using the open-source C++/Python finite element method package, \texttt{FEniCS} \citep{Logg2012}.
Adaptive time-stepping is used to constrain the local error to $\leq \num{e-4}$.
The second order Heun's method is used to determine the local error of the first order explicit method.
It was found in test cases that the largest source of local error is the velocity fields and that the local error between the velocity fields is comparable between tentative and updated velocities.
Therefore, only the tentative velocities are used in determining the local error and thus decreasing the number of second order solves to just one per time step.

\Cref{eq:continuity} is a pure advection equation, which is susceptible to node-to-node oscillations \citep{Brooks1982}.
To prevent the oscillations from occurring in the simulations, \cref{eq:alpha_update_discretized} is stabilized using the streamline-upwind/Petrov-Galerkin (SUPG) formulation for convection-dominated flows.
In the SUPG formulation, the test function in the finite element formulation is modified to allow for upwinding \citep{Brooks1982}:
\begin{equation}
\varphi' = \varphi + \tau_{SUPG} \vb*{v}_d \vdot \grad{\varphi}, \label{ch:results2:eq:SUPG_test_fn}
\end{equation}
where $\varphi$ is the test function and $\tau_{SUPG}$ is given by:
\begin{gather}
\tau_{SUPG} = \frac{h}{2\norm{\vb*{v}_d}} z, \label{ch:results2:eq:tau_SUPG} \\
z = \coth Pe_{\alpha} - \frac{1}{Pe_{\alpha}}. \label{ch:results2:eq:z_SUPG}
\end{gather}
\nomenclature{$Pe$}{P\'{e}clet number}%
\nomenclature{$h$}{Mesh element length}%
$h$ is the element length and $Pe$ is the P\'{e}clet number.
In pure advection transport, the P\'{e}clet number is infinite and $z$ is thus equal to one.


In order to maintain phase fraction boundedness, the IPCS scheme is bounded through the use of the nonlinear variational inequality solver SNES \citep{petsc-user-ref}.
In variational inequality, for a given lower and upper bounds, $l_i$ and $u_i$, and a continuously differentiable function, $F$, one of the following three cases will hold \citep{Munson2001}:
\begin{align}
x_i^* = l_i ~ &\text{then} ~ F_i(x^*) \geq 0,\\
x_i^* \in \qty(l_i, u_i) ~ &\text{then} ~ F_i(x^*) = 0,\\
x_i^* = u_i ~ &\text{then} ~ F_i(x^*) \leq 0,
\end{align}
where $x^*$ is the solution.
The algorithm used to solve this problem is the reduced-space method and is presented in \citet{Benson2006}.
\added{Using a variational form to impose inequality constraints is convenient when using the finite-element method in that, it is to implement the inequality constraints directly in the weak form and use a standard nonlinear solver.
The main drawback to this approach is the introduction of additional nonlinear constraints for each inequality constraint.}

An alternative method, a bound-constrained solver for linear variational inequality from the TAO suite \citep{tao-user-ref}, which solves the problem using the trust-region Newton method was also used for comparative purposes.
\added{This involves modifying the iterative linear step -- using an approximation of the Jacobian -- such that it will change the value of a bounded variable to be outside of the prescribed bounds.}
To ensure that the linear variation inequality solver converges at every time step, an additional constraint is placed on the time step such that the time step is decreased until the gradient of the objective function satisfies the specified tolerance.
However, this required the time step to consistently be in the order of $\num{e-8}$ while the SNES solver had no such requirement.
\added{This is likely due to the fact that convergence is dependent on approximation of both the Jacobian and Hessian matrices of the nonlinear system.}
Thus, the SNES solver is used in the bounded simulations presented in the subsequent sections.

\section{Results and Discussion}\label{sec:results}

In order to evaluate the effect of phase fraction boundedness on the IPCS scheme, simulations are initially performed with the assumption that the bulk and interfacial pressures are equal, $P_c = P_d = P_{int}$ (\cref{sec:TFM}).
These results are compared to each other and to an alternative finite-volume implementation of the two-fluid model in \texttt{OpenFOAM}, \texttt{twoPhaseEulerFoam}.
Following this, simulations are performed for the same conditions, but without the assumption that $P_c = P_{int}$, which is both a more accurate approximation and has been shown to increase the phase fraction interval over which the two-fluid model is well-posed \citep{Pauchon1986,Vaidheeswaran2016a,Stuhmiller1977}.

\subsection{Effects of Phase Fraction Boundedness, $P_c = P_{int}$} \label{sec:results1}

\begin{figure}
    \centering
    \begin{subfigure}{0.2\linewidth}
        \includegraphics[width=\linewidth]{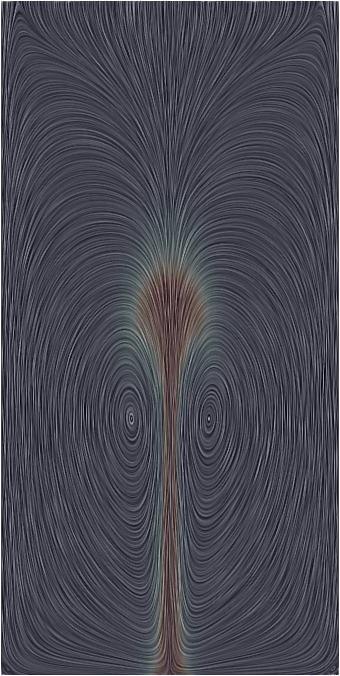}
        \caption*{1.25 s}
    \end{subfigure}
    \begin{subfigure}{0.2\linewidth}
        \includegraphics[width=\linewidth]{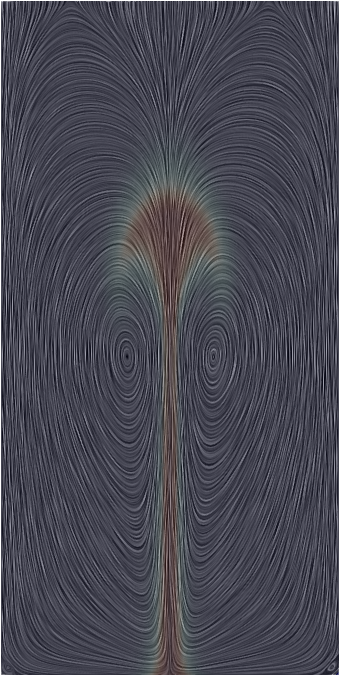}
        \caption*{1.41 s}
    \end{subfigure}
    \begin{subfigure}{0.2\linewidth}
        \includegraphics[width=\linewidth]{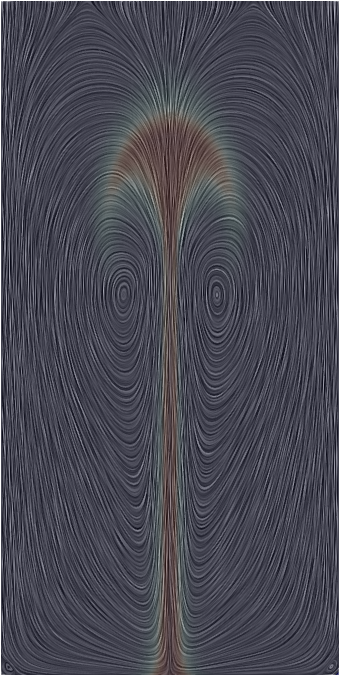}
        \caption*{1.56 s}
    \end{subfigure}
    \begin{subfigure}{0.2\linewidth}
        \includegraphics[width=\linewidth]{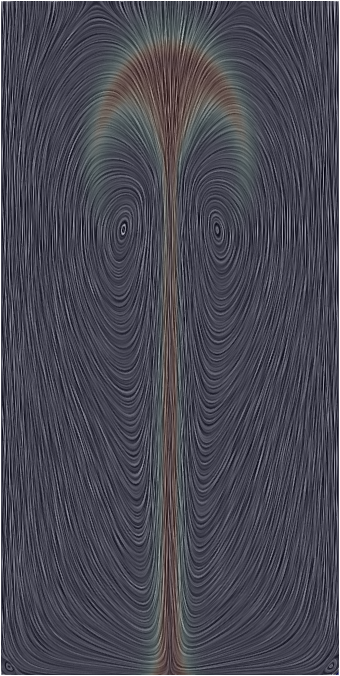}
        \caption*{1.72 s}
    \end{subfigure}
    \begin{subfigure}{0.076\linewidth}
        \includegraphics[width=\linewidth]{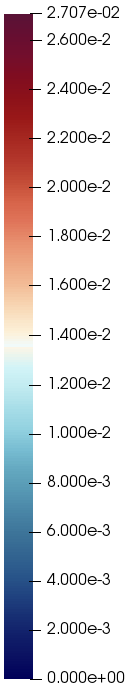}
        \caption*{}
    \end{subfigure}
    \caption{Evolution of the phase fraction and liquid velocity streamline over time. Colors denote $\alpha_g$.}
    \label{fig:var_CD}
\end{figure}

Simulations were performed using physical properties in \cref{tab:phys_prop} and auxiliary conditions in \cref{tab:bcs_var_CD} using the (i) unbounded IPCS solver, (ii) the bounded IPCS solver and (iii) the \texttt{twoPhaseEulerFoam} solver from \texttt{OpenFOAM}.
All of these solvers use the assumption that $P_c = P_{int}$.
\Cref{fig:var_CD} shows $\alpha_g$ and the liquid velocity streamlines at various times for the bounded simulation.
Simulation results are shown starting at $t=\SI{1.25}{\second}$, when the Rayleigh-Taylor instability \citep{Ruzicka2003} first manifests in the formation of a gas phase ``plume'' as it convects through the liquid phase.
As expected, as time increases the plume width increases with increasingly large vortices in the liquid phase velocity forming in its wake.
This morphology has been observed experimentally in the startup period of rectangular bubble columns \citep{Mudde2005} and its observation serves as qualitative experimental validation of the simulation results.
For all three simulations (unbounded IPCS, bounded IPCS, and \texttt{twoPhaseEulerFoam}), the bubble plume rises in the column center and vortices in the liquid velocity are observed on each side of the plume.

\begin{figure}
    \centering
    \begin{subfigure}{0.2\linewidth}
        \includegraphics[width=\linewidth]{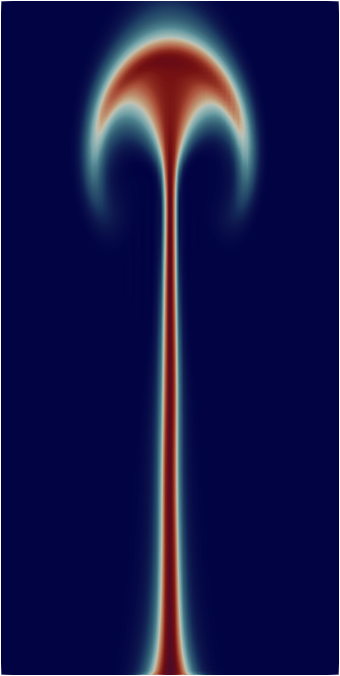}
        \caption*{}\label{fig:var_CD_unbounded}
    \end{subfigure}
    \begin{subfigure}{0.2\linewidth}
        \includegraphics[width=\linewidth]{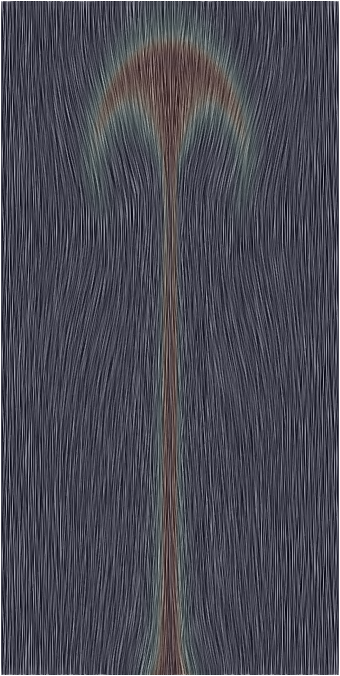}
        \caption*{}\label{fig:var_CD_unbounded_vg}
    \end{subfigure}
    \begin{subfigure}{0.2\linewidth}
        \includegraphics[width=\linewidth]{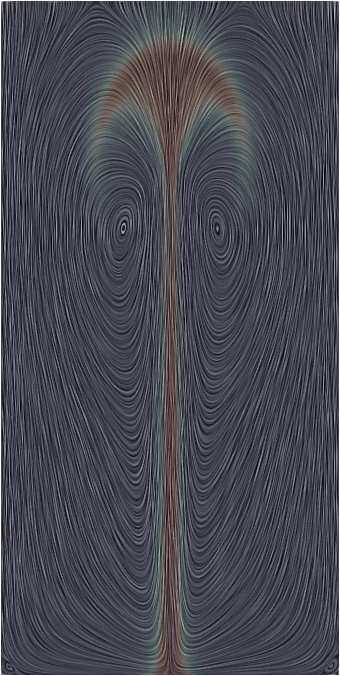}
        \caption*{}\label{fig:var_CD_unbounded_vl}
    \end{subfigure}
    \begin{subfigure}{0.077\linewidth}
        \includegraphics[width=\linewidth]{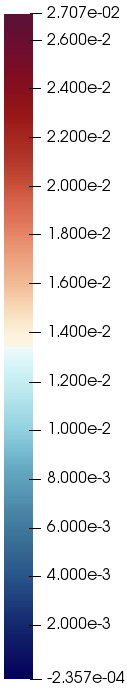}
        \caption*{}
    \end{subfigure}
    \\
    \begin{subfigure}{0.2\linewidth}
        \includegraphics[width=\linewidth]{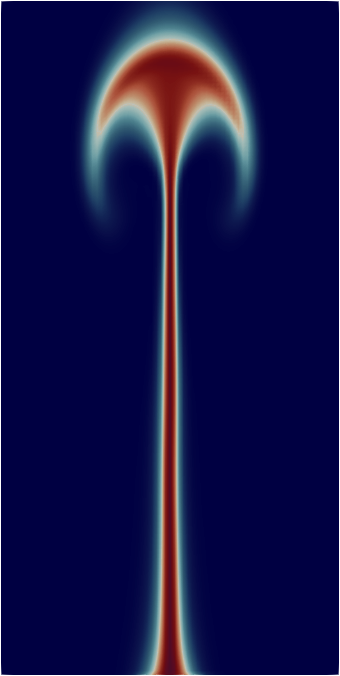}
        \caption*{}\label{fig:var_CD_bounded}
    \end{subfigure}
    \begin{subfigure}{0.2\linewidth}
        \includegraphics[width=\linewidth]{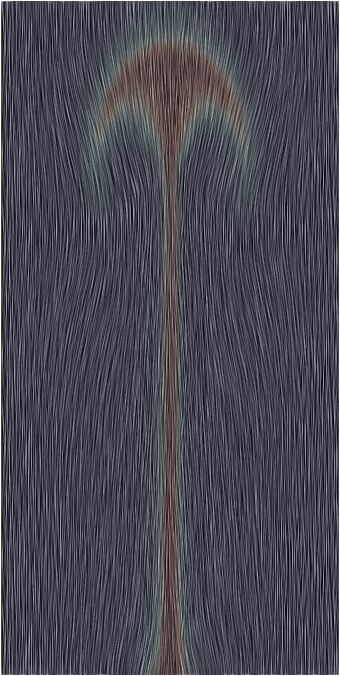}
        \caption*{}\label{fig:var_CD_bounded_vg}
    \end{subfigure}
    \begin{subfigure}{0.2\linewidth}
        \includegraphics[width=\linewidth]{var_CD/bounded_vl_550}
        \caption*{}\label{fig:var_CD_bounded_vl}
    \end{subfigure}
    \begin{subfigure}{0.076\linewidth}
        \includegraphics[width=\linewidth]{var_CD/bounded_cbar}
        \caption*{}
    \end{subfigure}
    \\
    \begin{subfigure}{0.2\linewidth}
        \includegraphics[width=\linewidth]{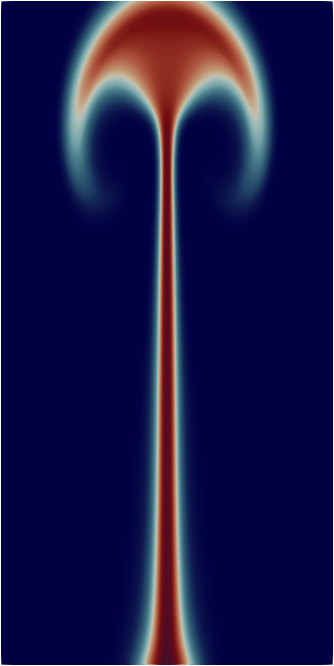}
        \caption*{$\alpha_g$}\label{fig:var_CD_OF}
    \end{subfigure}
    \begin{subfigure}{0.2\linewidth}
        \includegraphics[width=\linewidth]{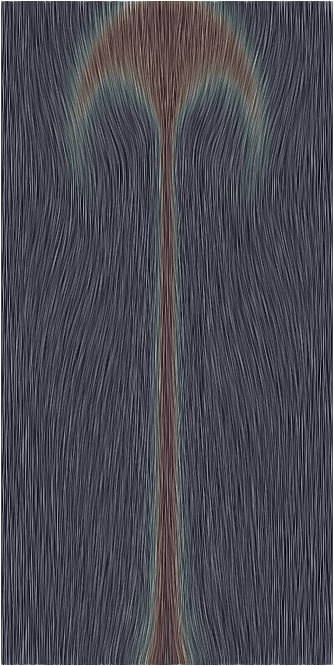}
        \caption*{$\vb*{v}_g$}\label{fig:OF_vg}
    \end{subfigure}
    \begin{subfigure}{0.2\linewidth}
        \includegraphics[width=\linewidth]{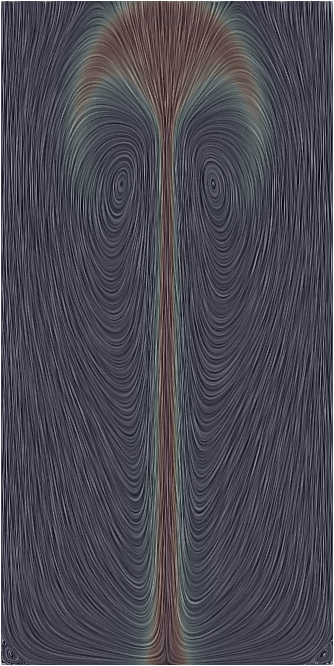}
        \caption*{$\vb*{v}_l$}\label{fig:OF_vl}
    \end{subfigure}
    \begin{subfigure}{0.0765\linewidth}
        \includegraphics[width=\linewidth]{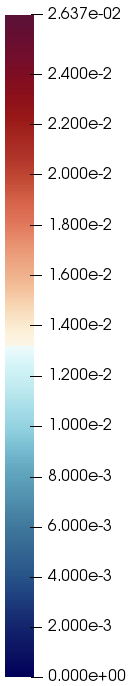}
        \caption*{}
    \end{subfigure}
    \caption{Surface plot of (left) phase fraction, (center) gas velocity and (right) liquid velocity at $t=\SI{1.72}{\second}$ from (top) unbounded IPCS, (middle) bounded IPCS and (bottom) \texttt{twoPhaseEulerFoam}.}
    \label{fig:var_CD_all}
\end{figure}

The average bubble Reynolds numbers at $t=\SI{1.72}{\second}$ computed using the definition given in \cref{tab:phys_prop} for the unbounded IPCS, bounded IPCS and \texttt{twoPhaseEulerFoam} are \numlist{11.542;11.369;11.724}, respectively.
The gas and liquid volumetric fluxes ($\alpha_g v_g$ and $\alpha_l v_l$) from the simulations ranges from \SIrange{4.93e-3}{5.11e-3}{\meter/\second} and \SIrange{0.126}{0.142}{\meter/\second}, respectively.
These values are well below the range of fluxes where bubble-induced turbulence was observed by past experimental studies \cite{Rzehak2013,Vaidheeswaran2017b}, supporting the choice to neglect the contribution of bubble-induced turbulence in this work.
The results from unbounded IPCS, bounded IPCS and \texttt{twoPhaseEulerFoam} solvers at $t=\SI{1.72}{\second}$ are given in \cref{fig:var_CD_all}.
From \cref{fig:var_CD_all}, both the unbounded and bounded IPCS solvers resulted in similar flow profiles and a qualitatively similar plume is also observed in the \texttt{twoPhaseEulerFoam} simulation.
However, the gas plume appears to be rising at a faster rate and is wider in the \texttt{twoPhaseEulerFoam} simulation than in the IPCS simulations.

In order to partially validate the simulations, an empirical correlation determined from multiple experimental studies for terminal velocity of a rigid sphere in a fluid \cite{Clift1978}, which predicts the terminal Reynolds number ($Re_T = \rho_l d v_T/\mu_l$) was used:
\begin{equation}
\log Re_T = -1.7095 + 1.33438\log N_D - 0.11591\qty(\log N_D)^2,
\end{equation}
where $N_D = 4\rho_l\Delta\rho g d^3/3\mu_l^2$.
A terminal Reynolds number of $Re_T = 11.441$ was obtained for the physical properties used in this study.
The terminal velocity resulting from the empirical correlation is $v_T=\SI{0.0572}{\meter/\second}$.
The average slip velocity computed from simulations is $v_T = \SI{0.0589}{\meter/\second}$ for the IPCS method and $v_T = \SI{0.0592}{\meter/\second}$ for the \texttt{OpenFOAM} solver.
This supports the numerical accuracy of both simulations, where the slip velocity should be similar to the terminal velocity predicted by the empirical correlation consistent with the parameters used in the drag force (Section \ref{sec:drag_force}).

\begin{figure}
    \centering
    \begin{subfigure}{0.45\linewidth}
        \includegraphics[width=\linewidth]{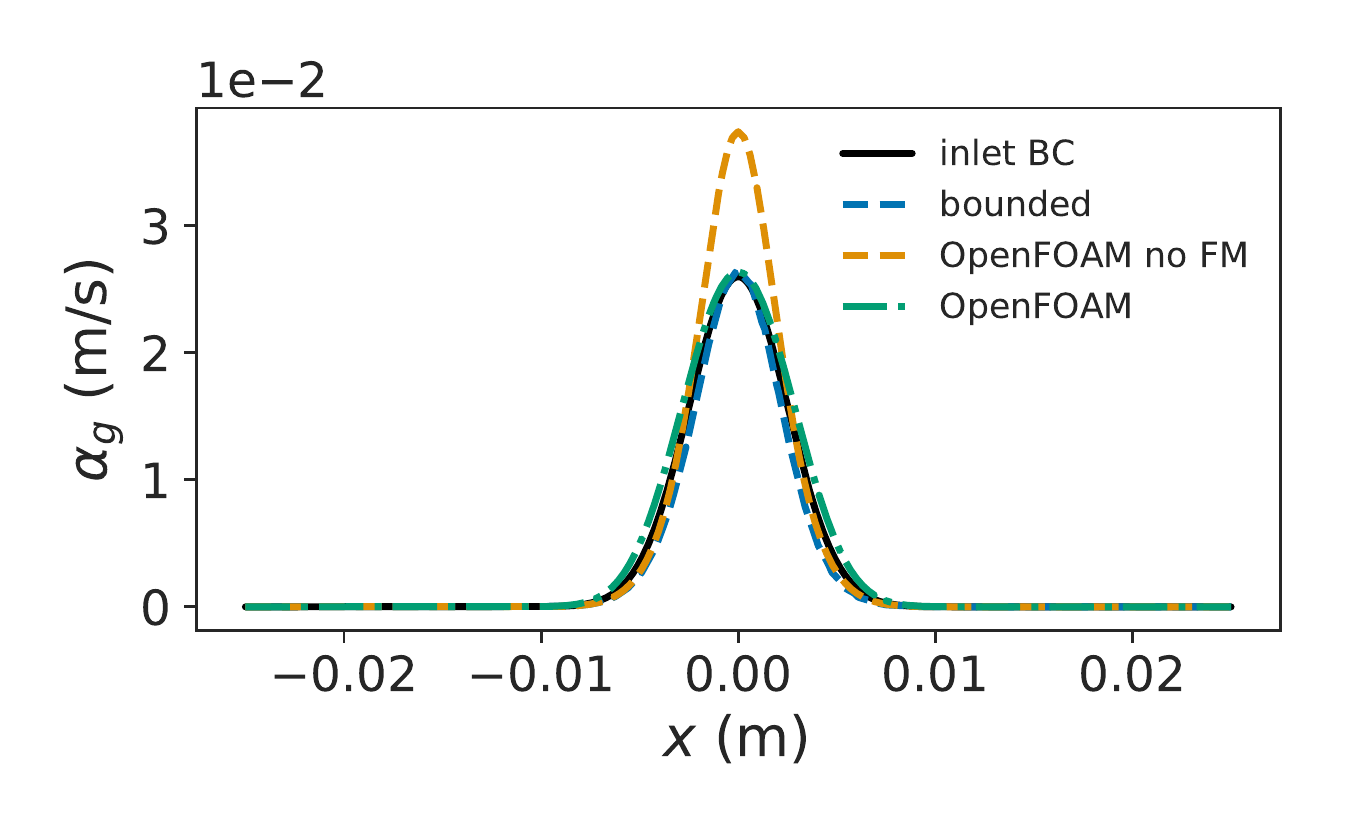}
        \caption{}
        \label{fig:velocity_diff_near_inlet_alpha_g}
    \end{subfigure}\\
    \begin{subfigure}{0.45\linewidth}
        \includegraphics[width=\linewidth]{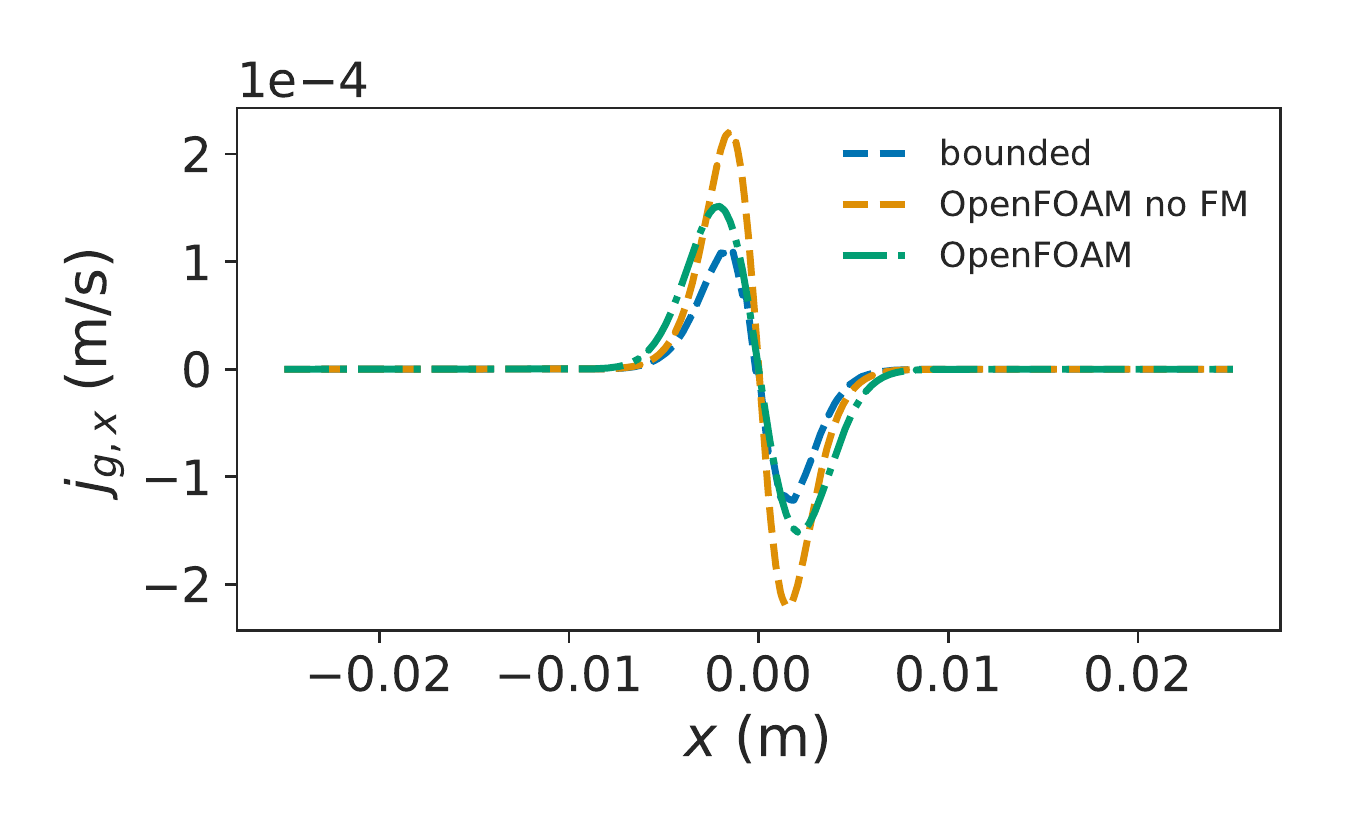}
        \caption{}
        \label{fig:velocity_diff_near_inlet_jx}
    \end{subfigure}
    \begin{subfigure}{0.45\linewidth}
        \includegraphics[width=\linewidth]{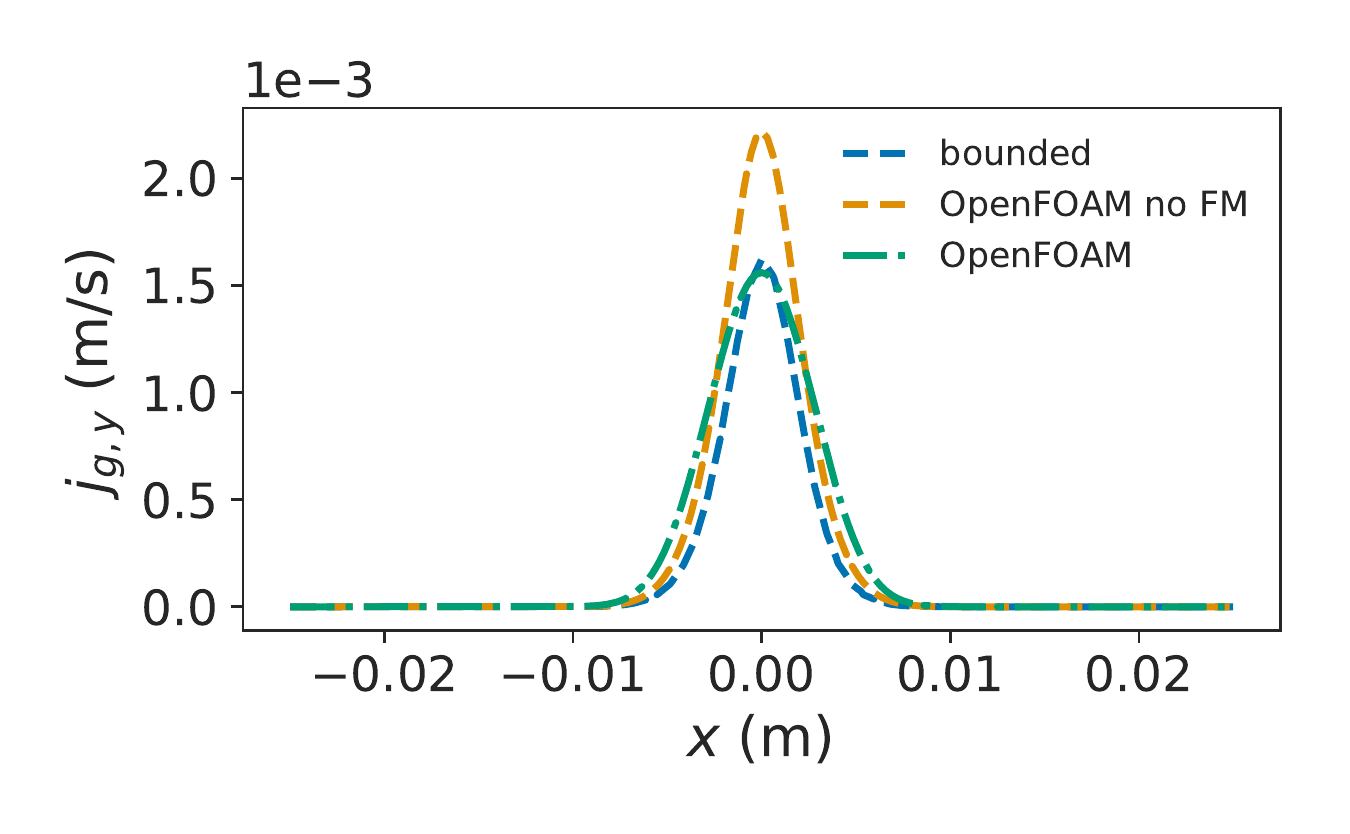}
        \caption{}
        \label{fig:velocity_diff_near_inlet_jy}
    \end{subfigure}
    \caption{Plots of the (a) gas fraction, (b) $x$- and (c) $y$-components of the gas volumetric flux, $\vb*{j}_g$, along $x$ at $y=\SI{0.0005}{\meter}$. FM denotes face-based momentum formulation.}
    \label{fig:velocity_diff_near_inlet}
\end{figure}

\Cref{fig:velocity_diff_near_inlet} shows the gas volumetric flux, $\vb*{j}_g$, and the gas phase fraction profile near the inlet region along $x$ for the bounded IPCS result and \texttt{OpenFOAM} using both the cell-based and face-based momentum formulations.
The cell-based formulation is the default formulation in \texttt{twoPhaseEulerFoam} and has been reported to exhibit numerical artifacts. 
In this work, the gas fraction and the gas volumetric flux are affected, where a distinct jump in the phase fraction profile from the inlet to inside the domain is observed in \cref{fig:velocity_diff_near_inlet_alpha_g}.
The \texttt{twoPhaseEulerFoam} results shown in \cref{fig:var_CD_all} are from the face-based formulation, where improved numerical stability is achieved by interpolating the momentum transport terms to the faces of the cells.

From \cref{fig:velocity_diff_near_inlet}, the change in the gas fraction profile as the gas enters the domain is significantly more gradual than the cell-based formulation and the profiles are in better agreement with the IPCS results.
However, it has been shown in \citet{Agnaou2018} that upstream flow conditions affect the flow profile downstream and using different inlet conditions will result in markedly different flow profiles, even at steady-state and with a sufficiently long pipe/channel.
It then stands to reason that difference between the gas volumetric flux profiles of the \texttt{twoPhaseEulerFoam} and IPCS simulations (\cref{fig:velocity_diff_near_inlet_jx,fig:velocity_diff_near_inlet_jy}) will result in a bubble plume that is wider and faster for \texttt{twoPhaseEulerFoam}.

\begin{figure}
    \centering
    \includegraphics[width=0.75\textwidth]{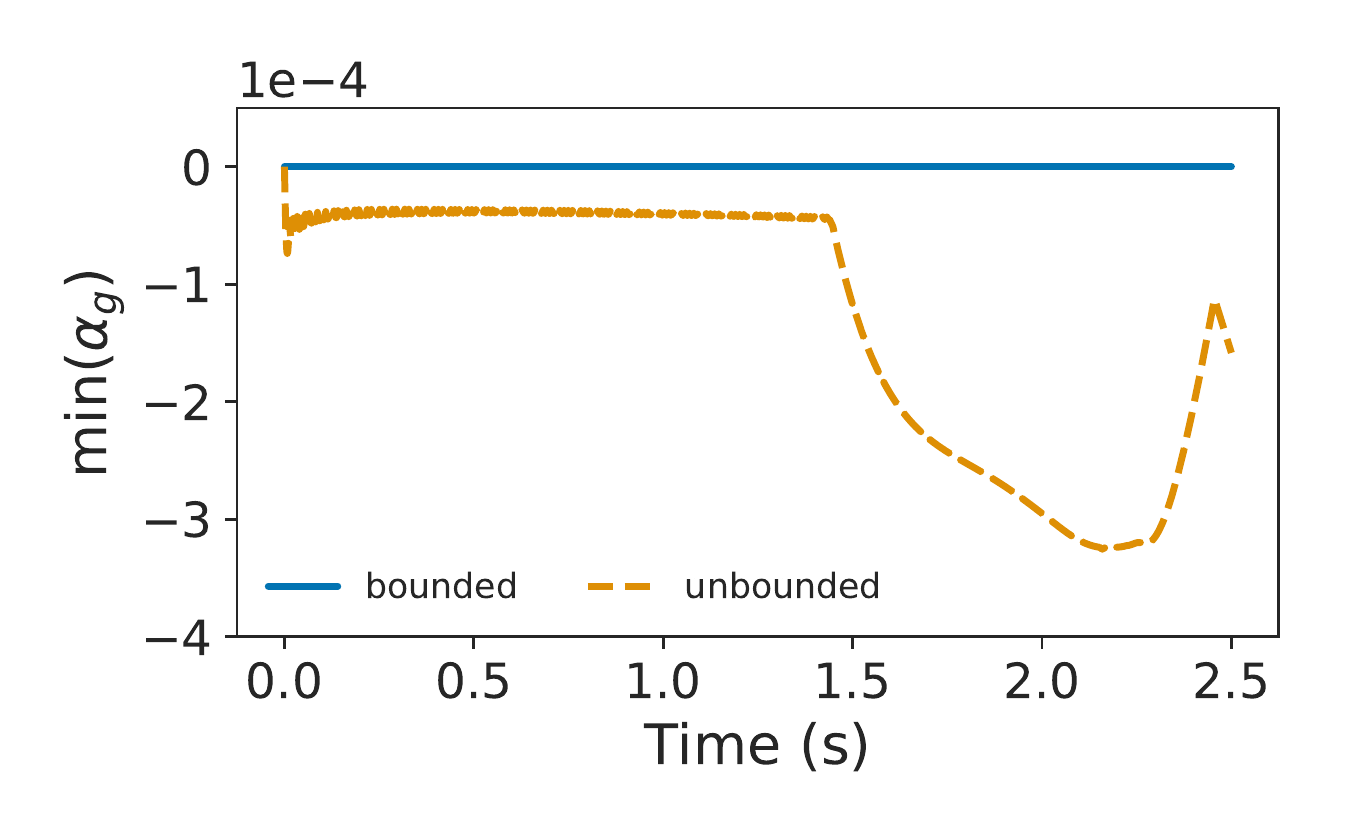}
    \caption{Evolution of $\min(\alpha_g)$ over time}
    \label{fig:var_min_alpha_g}
\end{figure}

\Cref{fig:var_min_alpha_g} shows the evolution of the minimum gas phase fraction, $\min{(\alpha_g)}$, for the unbounded and bounded IPCS solvers.
For the unbounded solver, the magnitude of $\min(\alpha_g)$ is found to be on the order of $\num{e-4}$, which is on the order of the relative error tolerance of the adaptive time-stepping method, but well above that of the underlying linear solver ($\num{e-13}$).
For the bounded solver, the magnitude is within the tolerance of the nonlinear variational solver ($\num{e-11}$).

\begin{figure}
    \centering
    \begin{subfigure}{0.48\linewidth}
        \includegraphics[width=\linewidth]{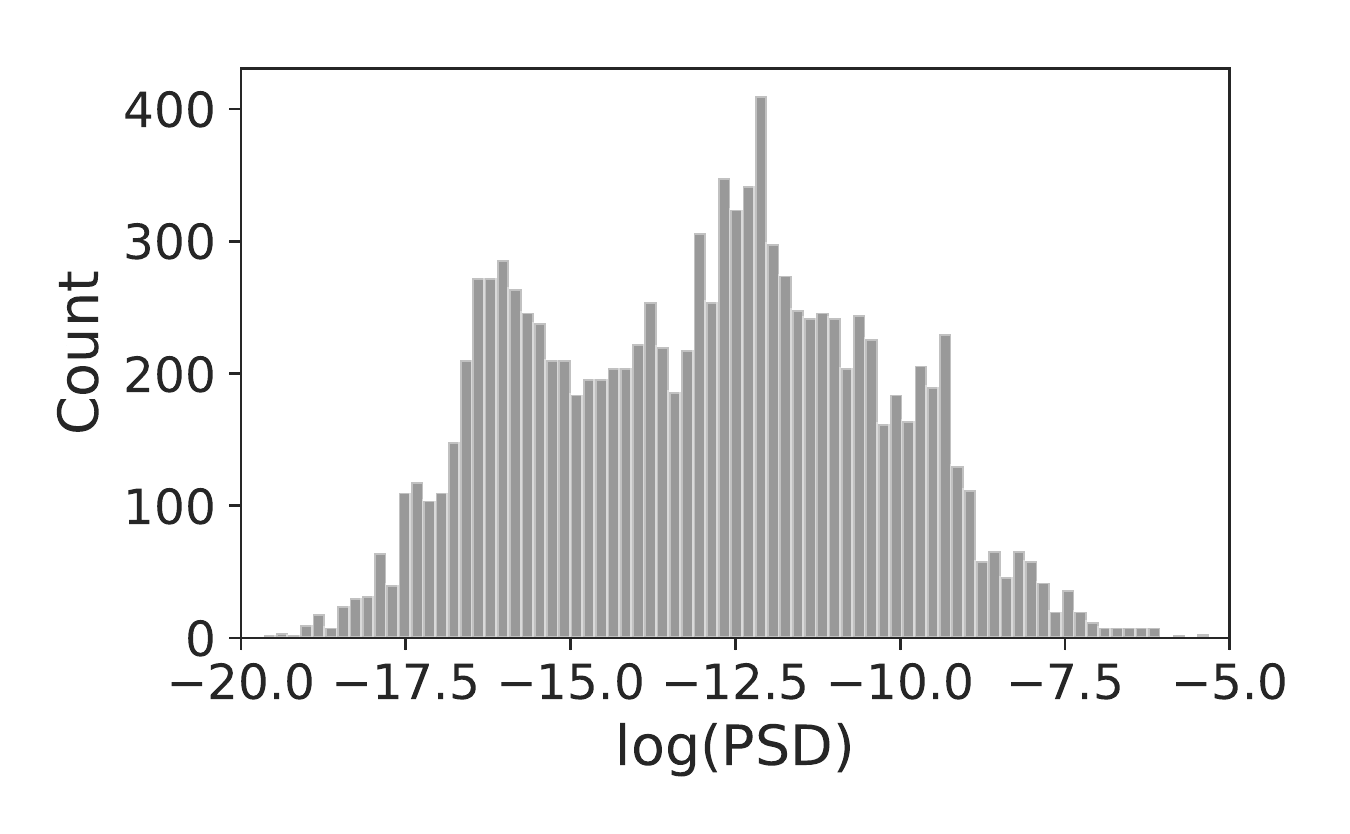}
        \caption{}\label{fig:fft_unbounded}
    \end{subfigure}
    \begin{subfigure}{0.48\linewidth}
        \includegraphics[width=\linewidth]{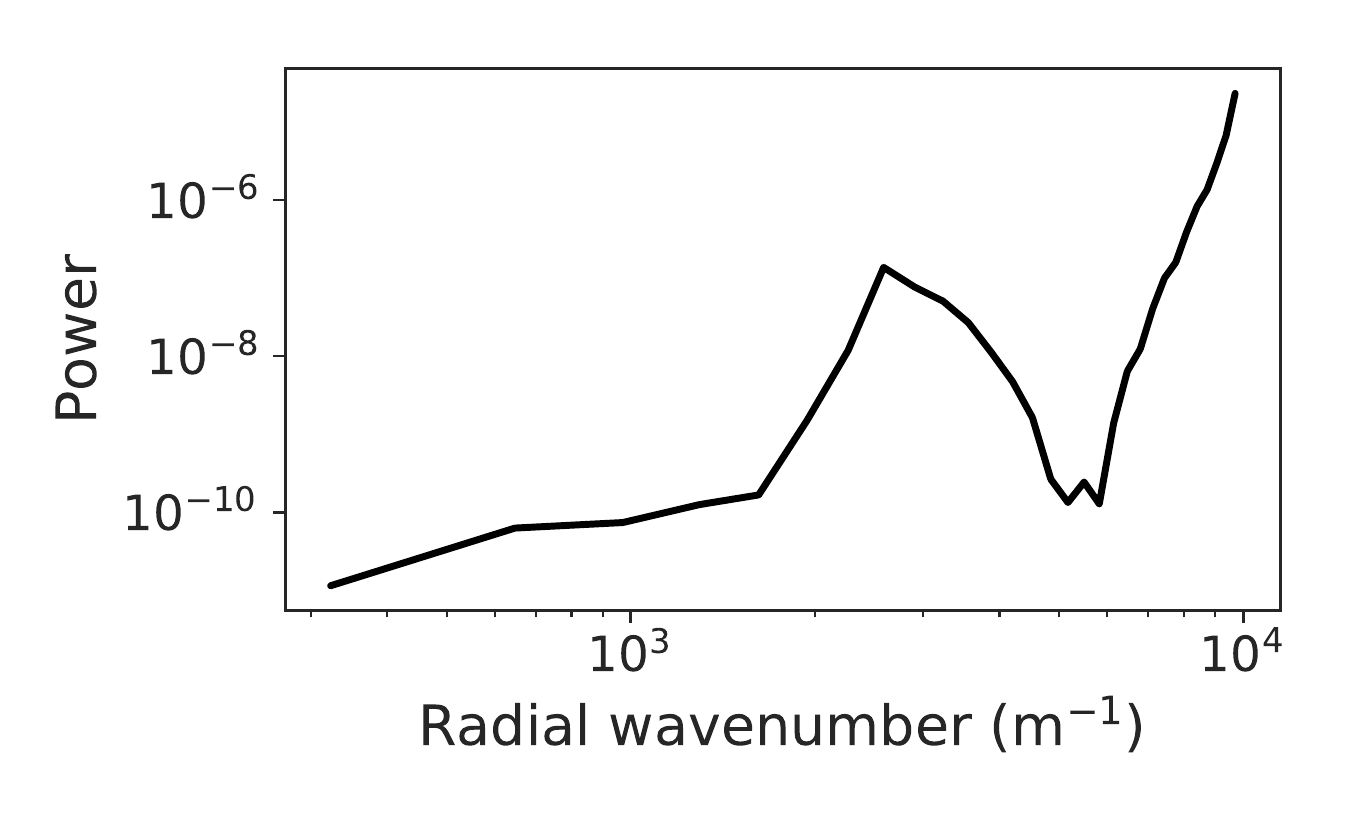}
        \caption{}\label{fig:fft_unbounded_radial}
    \end{subfigure}
    \\
    \begin{subfigure}{0.48\linewidth}
        \includegraphics[width=\linewidth]{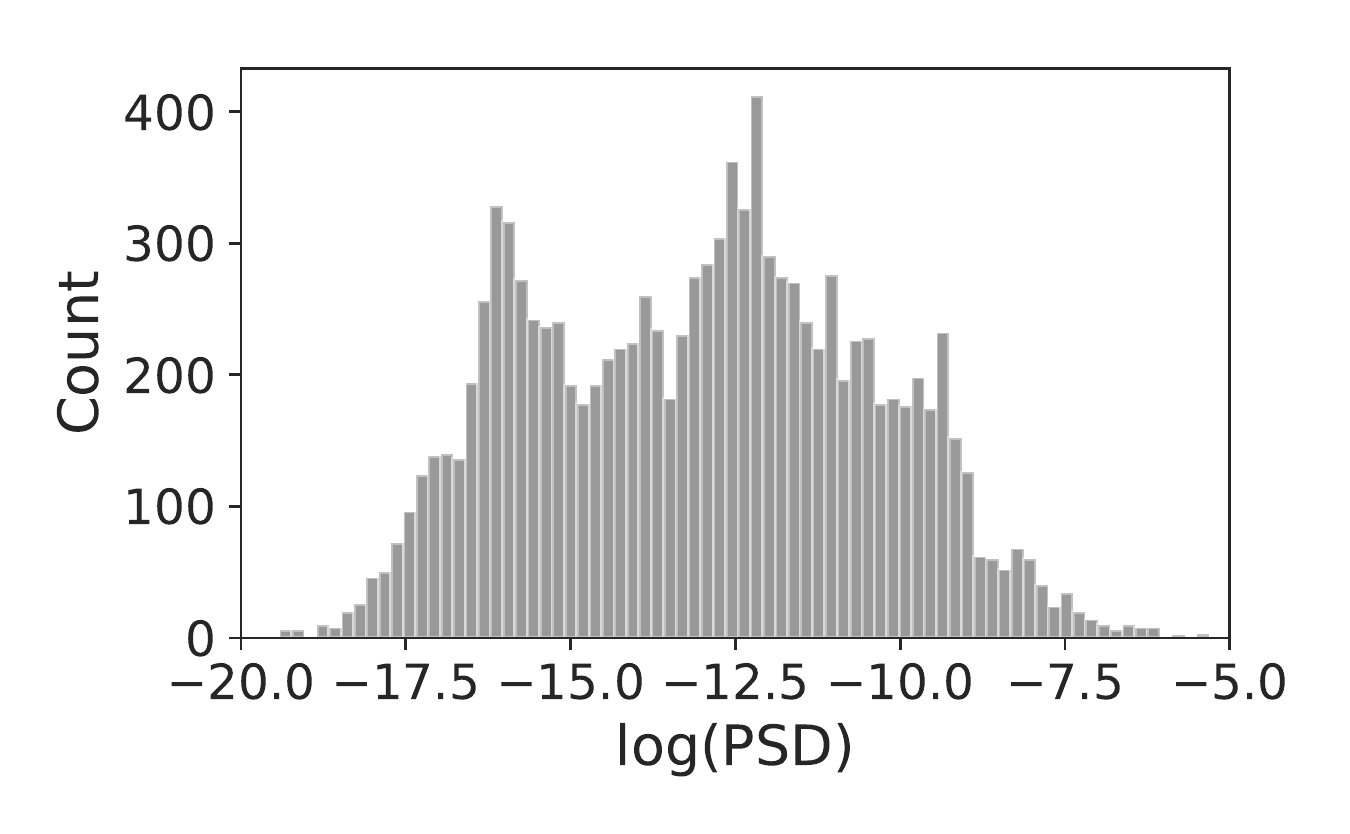}
        \caption{}\label{fig:fft_bounded}
    \end{subfigure}
    \begin{subfigure}{0.48\linewidth}
        \includegraphics[width=\linewidth]{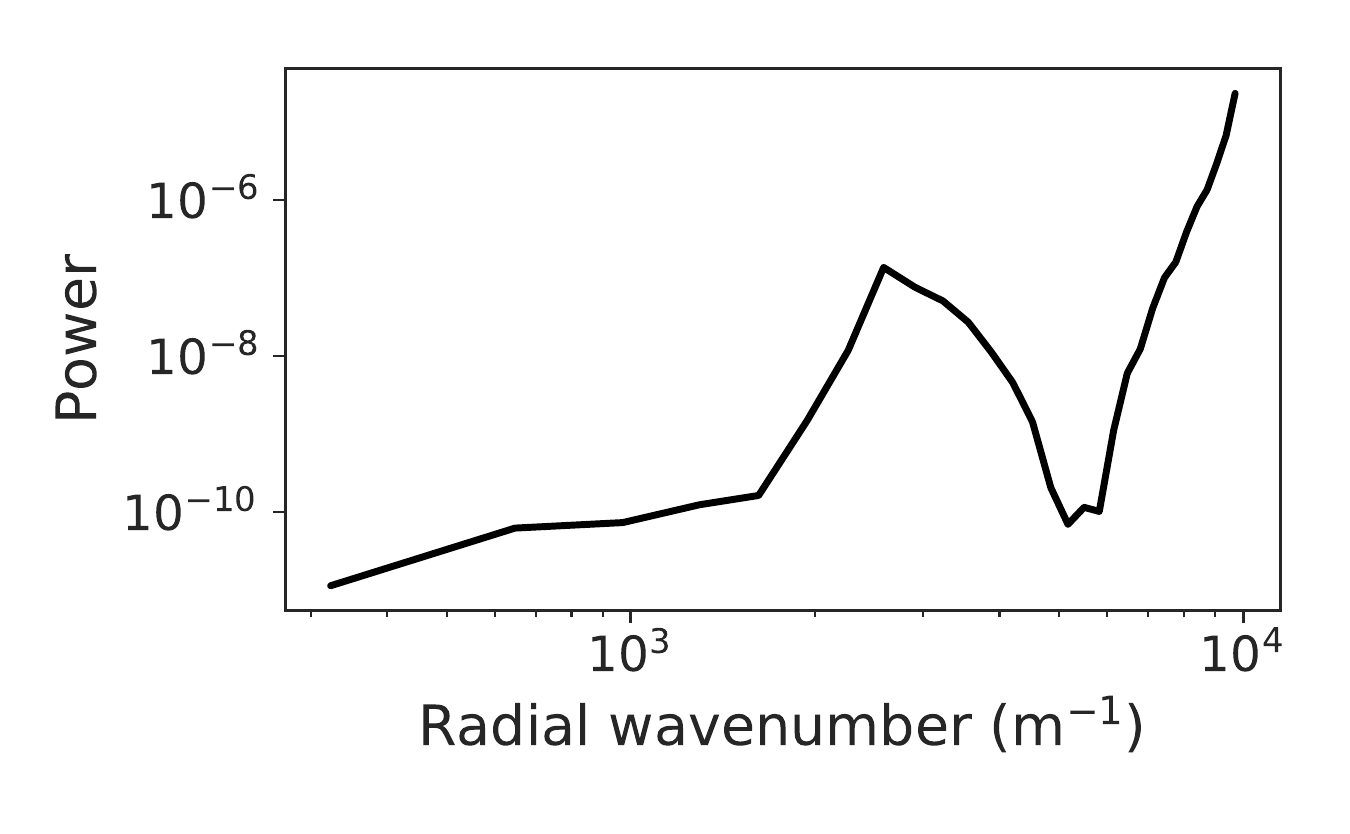}
        \caption{}\label{fig:fft_bounded_radial}
    \end{subfigure}
    \\
    \begin{subfigure}{0.48\linewidth}
        \includegraphics[width=\linewidth]{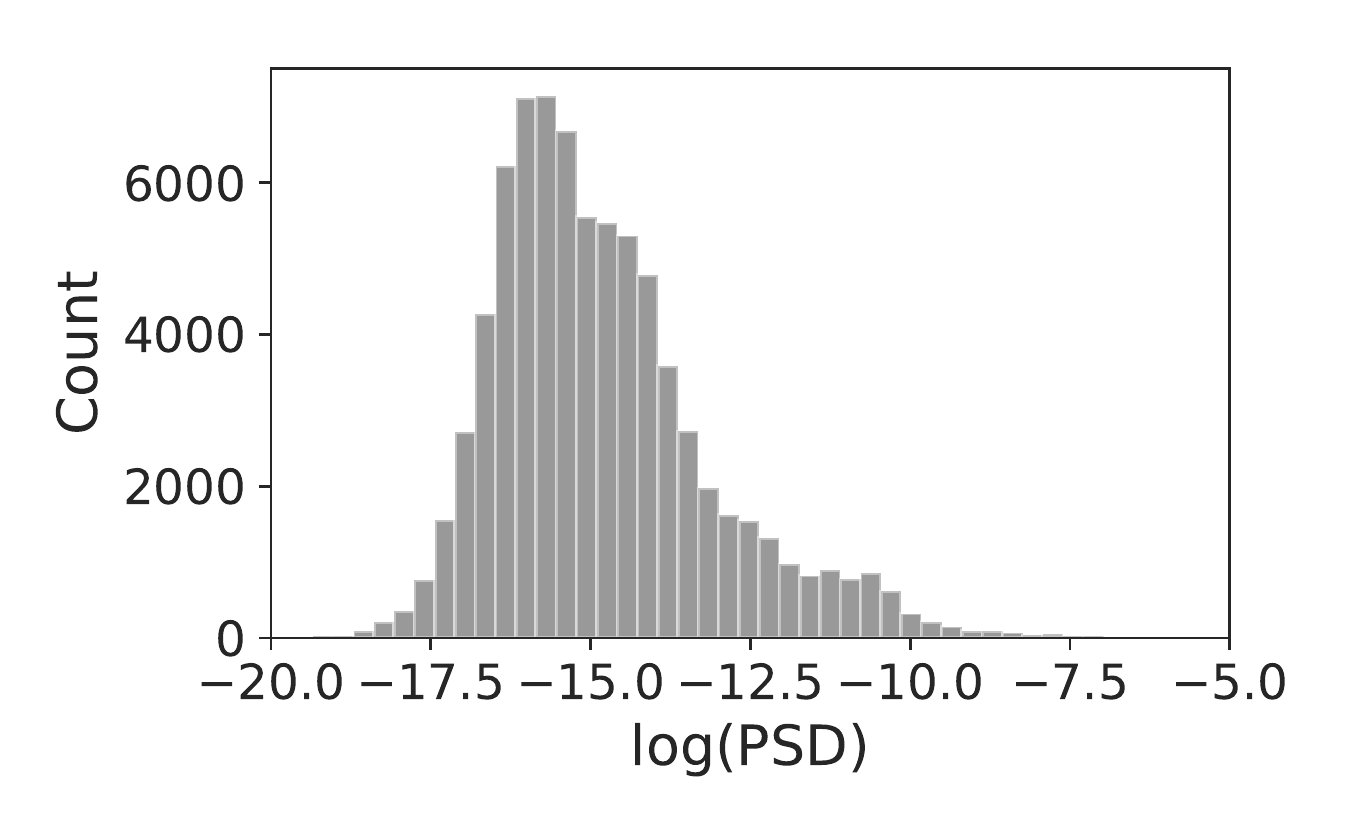}
        \caption{}\label{fig:fft_OF}
    \end{subfigure}
    \begin{subfigure}{0.48\linewidth}
        \includegraphics[width=\linewidth]{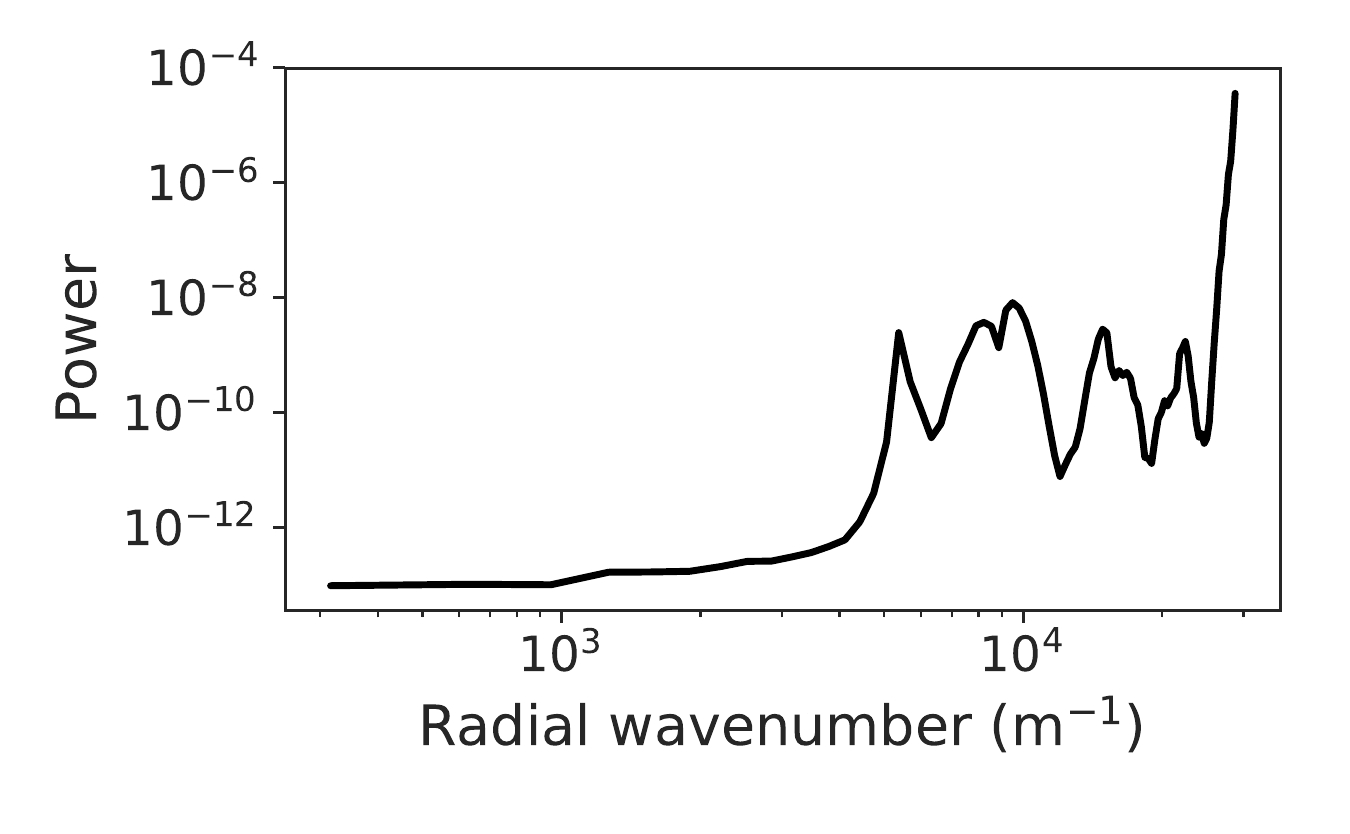}
        \caption{}\label{fig:fft_OF_radial}
    \end{subfigure}
    \caption{Two-dimensional spectral analysis of the phase fraction at $t=\SI{1.72}{\second}$ from (top) unbounded IPCS, (middle) bounded IPCS and (bottom) \texttt{twoPhaseEulerFoam} solvers. Left: histogram of power spectral density. Right: radially-averaged power spectrum}
    \label{fig:fft}
\end{figure}

Quantitative comparison of structure of the simulation results from the three solvers was performed using spectral analysis of the phase fraction profiles.
\Cref{fig:fft} contains the results from spectral analysis at $t=\SI{1.72}{\second}$ for the three different solvers.
The histograms of the power spectral density of the phase fraction computed using the IPCS solvers also contained near-zero (in the order of $\num{e-40}$) power densities that were omitted from \cref{fig:fft_unbounded,fig:fft_bounded}.
The power spectral density distributions of the IPCS results are similar but drastically different from that of the \texttt{twoPhaseEulerFoam} results.
This is corroborated by the power spectra shown in \cref{fig:fft_unbounded_radial,fig:fft_bounded_radial,fig:fft_OF_radial}.
Again, the deviations are possibly due to both the phase fraction different approaches to phase fraction boundedness and/or the spatial interpolation schemes.

\subsection{Effects of Phase Fraction Boundedness $P_c \ne P_{int}$} \label{sec:results2}

\begin{figure}
    \centering
    \begin{subfigure}{0.2\linewidth}
        \includegraphics[width=\linewidth]{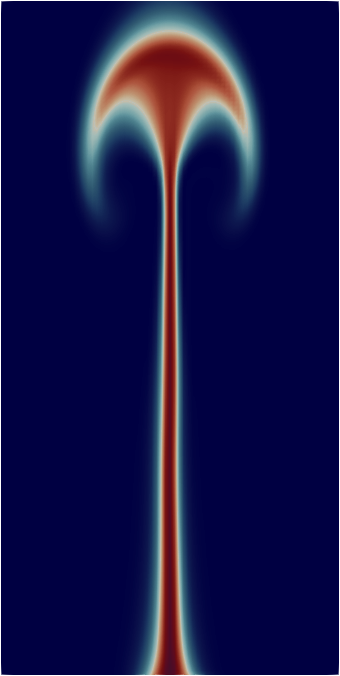}
        \caption*{$\alpha_g$}\label{fig:var_CD_bounded_intp}
    \end{subfigure}
    \begin{subfigure}{0.2\linewidth}
        \includegraphics[width=\linewidth]{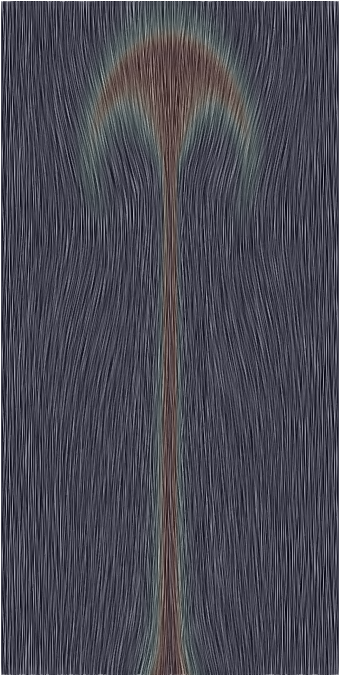}
        \caption*{$\vb*{v}_g$}\label{fig:var_CD_bounded_intp_vg}
    \end{subfigure}
    \begin{subfigure}{0.2\linewidth}
        \includegraphics[width=\linewidth]{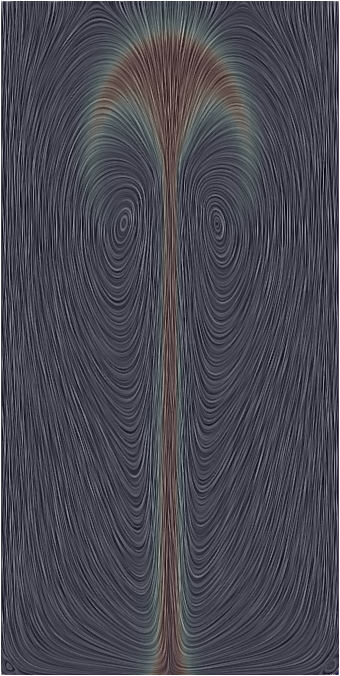}
        \caption*{$\vb*{v}_l$}\label{fig:var_CD_bounded_intp_vl}
    \end{subfigure}
    \begin{subfigure}{0.076\linewidth}
        \includegraphics[width=\linewidth]{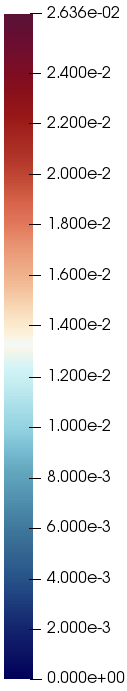}
        \caption*{}
    \end{subfigure}
    \caption{Surface plot of (left) phase fraction, (center) gas velocity and (right) liquid velocity at $t=\SI{1.72}{\second}$ from bounded IPCS with interfacial pressure.}
\label{fig:var_CD_intp}
\end{figure}

\begin{figure}
    \centering
    \begin{subfigure}{0.48\linewidth}
        \includegraphics[width=\linewidth]{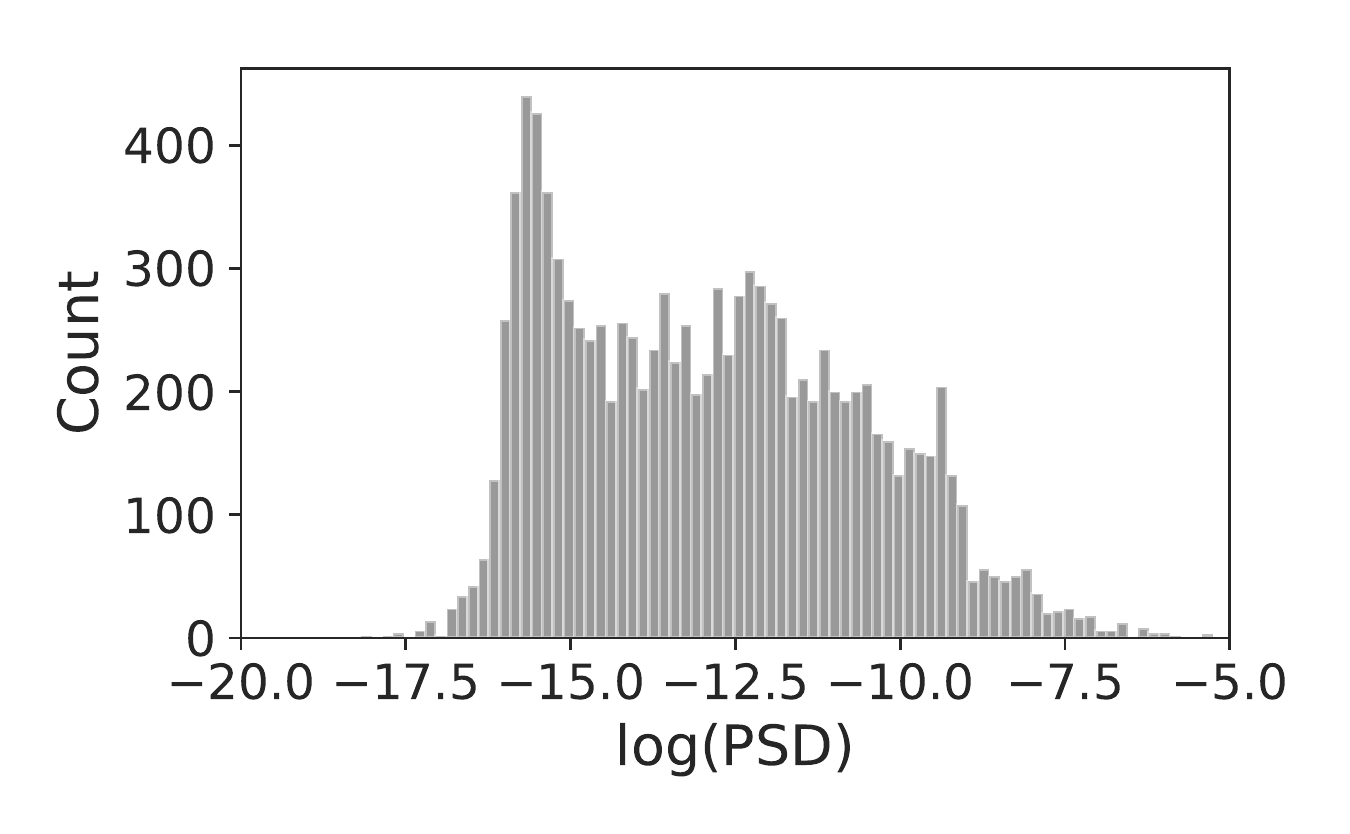}
        \caption{}\label{fig:fft_bounded_intp}
    \end{subfigure}
    \begin{subfigure}{0.48\linewidth}
        \includegraphics[width=\linewidth]{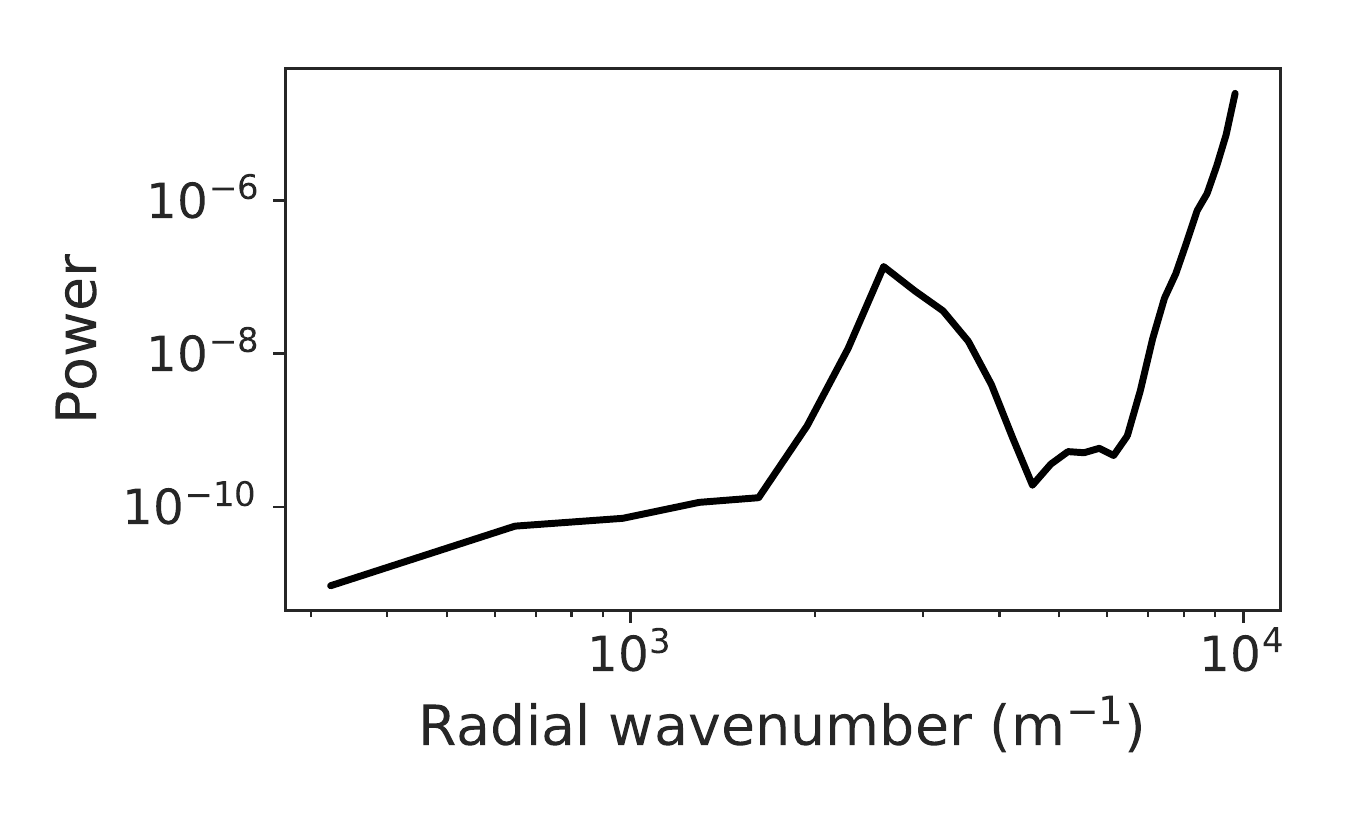}
        \caption{}\label{fig:fft_bounded_intp_radial}
    \end{subfigure}
    \caption{Two-dimensional spectral analysis of the phase fraction at $t=\SI{1.72}{\second}$ from bounded IPCS with interfacial pressure. Left: histogram of power spectral density. Right: radially-averaged power spectrum}
\label{fig:fft2}
\end{figure}

The previous simulations were performed with the assumption that the interfacial pressure and the bulk pressure of the continuous phase are equal to each other.
In this section, this assumption is removed and a bounded IPCS simulation is performed to assess its effect, which has been shown to increase the phase fraction interval over which the two-fluid model is well-posed \citep{Pauchon1986,Vaidheeswaran2016a,Stuhmiller1977}.
The flow profile at $t=\SI{1.72}{\second}$ is shown in \cref{fig:var_CD_intp} and is qualitatively similar to the results in \cref{fig:var_CD_all}.
The spatial variation of the gas fractions within the gas phase plume appears to be smoother than the results from the bounded IPCS solver in the previous section (\cref{fig:var_CD_all}).
The gas and liquid volumetric fluxes are \SI{4.79e-3}{\meter/\second} and \SI{0.130}{\meter/\second}, respectively.
The volumetric fluxes fall inside the range of the reported values from simulations performed in \cref{sec:results1}.
Additionally, the terminal velocity determined from the simulation is \SI{0.0592}{\meter/\second}, the same as previous values from \cref{sec:results1}.
The histogram of the power spectral density (\cref{fig:fft_bounded_intp}) also shows qualitative agreement with \cref{fig:fft_bounded} but the radially-averaged power spectra of the two cases are slightly different (\cref{fig:fft_bounded_radial,fig:fft_bounded_intp_radial}).

\begin{figure}
    \centering
    \includegraphics[width=0.75\linewidth]{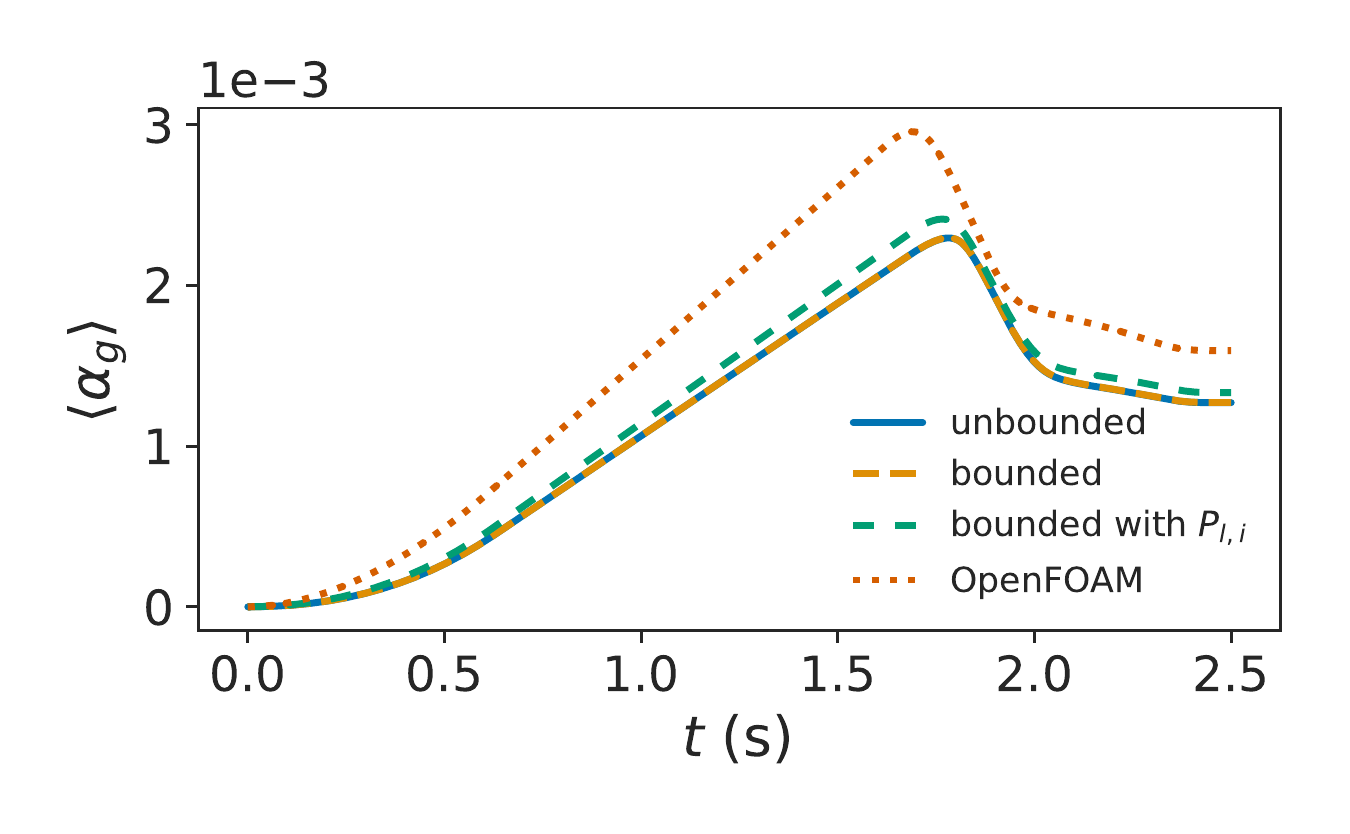}
    \caption{Time-evolution of gas holdup with variable $C_D$}\label{fig:holdup_varCD}
\end{figure}

The time evolution of the overall gas holdup, $\langle\alpha_g\rangle$, from the simulations with the assumption $P_c = P_{int}$ (\cref{sec:results1}) and the simulation with $P_c \neq P_{int}$ are compared in \cref{fig:holdup_varCD}.
The overall holdup from the IPCS solvers show little variation from each other, with some expected variation between the solutions with $P_c = P_{int}$ and $P_c \ne P_{int}$.
However, the gas holdup obtained from \texttt{twoPhaseEulerFoam} is consistently higher, although evolves in a qualitatively similar manner as the IPCS variants.
The peak observed in \cref{fig:holdup_varCD} corresponds to the point where the bubble plume is the largest, which is also right before the plume starts to exit the simulation domain.
Compared to the IPCS simulations, the time in which the peak occurs is earlier for the \texttt{twoPhaseEulerFoam} simulation.
This supports the qualitative observations made that the size of the plume and the rate in which the plume moves through the liquid are similar among the IPCS simulations but different when simulations under the same conditions that are performed using \texttt{twoPhaseEulerFoam}.

\section{Conclusions} \label{sec:conclusions}

A phase-bounded numerical method for the two-fluid model is developed using the incremental pressure correction scheme.
The phase fraction boundedness is imposed implicitly through the use of the SNES variational inequality solver.
Simulations are performed to compare the solution obtained from the phase-bounded method and are found to be similar to the unbounded method, but with the phase fraction equality constraints satisfied within the tolerance of the nonlinear variational inequality solver.
The results from the unbounded method exhibit deviations of the minimum value of the gas phase fraction in the domain that are several orders of magnitude greater than the linear solver error tolerance.

Qualitative agreement of the flow profile is found with respect to an alternative bounded two-fluid model solver, the \texttt{twoPhaseEulerFoam} solver in \texttt{OpenFOAM}, although quantitative agreement is limited to the slip velocity.
This is attributed to either or both the difference in method for imposing phase fraction bounds and approach to spatial interpolation.
All numerical solutions are found to agree qualitatively with experimental studies of two-dimensional rectangular bubble columns in the literature.

Finally, the effects of the assumption of the interfacial and bulk pressures of the continuous phase being equal were studied using the phase-bounded method and found to be non-negligible.

\section*{Acknowledgments}

This research was supported by the Natural Sciences and Engineering Research Council (NSERC) of Canada and Compute Canada.

\bibliographystyle{elsarticle-harv}
\bibliography{multiphase,computational}

\newpage
\begin{appendices}
\nomenclature{$\varphi$}{Test function}%
\section{Weak Form of Governing Equations}\label{sec:weak_forms}
The time-discretized equations presented in the previous section are solved using the method of lines.
The finite element method is used to approximate the spatial derivatives in the differential equations.
The finite-element method requires the differential equations to be formulated into their weak formulations.
In this section, weak formulations of \cref{eq:mom_g_discretized,eq:mom_l_discretized,eq:pp_discretized,eq:vl_update_discretized,eq:vg_update_discretized,eq:alpha_update_discretized} with the appropriate boundary conditions are described.
The time discretization scheme is still the explicit Euler method but extension to other time discretization schemes is possible.

\paragraph{Tentative Velocity}
Taking the inner product of \cref{eq:mom_g_discretized,eq:mom_l_discretized} with the test function for each phase, $\vb*{\varphi}_q$, and integrating over the domain yields the following:
\begin{subequations}
    \begin{align}
    \begin{split}
    \inp*{\frac{\tilde{\vb*{v}}_{l}^{*}- \tilde{\vb*{v}}_{l}^n}{\Delta t}}{\vb*{\varphi}_l}_{\Omega} + \inp*{\tilde{\vb*{v}}_{l}^n\vdot\tilde{\grad} \tilde{\vb*{v}}_{l}^n}{\vb*{\varphi}_l}_{\Omega} &= -\inp*{Eu_l \tilde{\grad} \tilde{P}_l^n}{\vb*{\varphi}_l}_{\Omega}  + \inp*{\frac{1}{Re_l}\frac{\tilde{\grad}\alpha_l^n \vdot \tilde{\vb*{\tau}}_l^{n+\frac{1}{2}}}{\alpha_l^n}}{\vb*{\varphi}_l}_{\Omega} \\
    & \quad - \inp*{\frac{1}{Re_l}\tilde{\vb*{\tau}}_l^{n+\frac{1}{2}}}{\grad \vb*{\varphi}_l}_{\Omega} + \inp*{\frac{1}{Re_l} \vb*{n} \vdot \tilde{ \vb*{\tau}}_l^{n+\frac{1}{2}}}{\vb*{\varphi}_l}_{\Gamma_N}\\
    & \quad + \inp*{\frac{1}{Fr^2}  \tilde{\vb*{g}}}{\vb*{\varphi}_l}_{\Omega}  +\inp*{\frac{3}{4}\frac{\alpha_g^n}{\alpha_l^n}\frac{C_{D}}{\tilde{d}_b}\norm{\tilde{\vb*{v}}_r^n}\tilde{\vb*{v}}_r^n}{\vb*{\varphi}_l}_{\Omega} \\
    & \quad -\inp*{ C_P \tilde{\vb*{v}}_r^n \vdot\tilde{ \vb*{v}}_r^n\frac{\tilde{\grad}\alpha_l^n}{\alpha_l^n}}{\vb*{\varphi}_l}_{\Omega},
    \end{split}\\
    \begin{split}
    \inp*{\frac{\tilde{\vb*{v}}_{g}^{*}- \tilde{\vb*{v}}_{g}^n}{\Delta t}}{\vb*{\varphi}_g}_{\Omega} + \inp*{\tilde{\vb*{v}}_{g}^n\vdot\tilde{\grad} \tilde{\vb*{v}}_{g}^n}{\vb*{\varphi}_g}_{\Omega} &= -\inp*{Eu_g \tilde{\grad}\qty(\tilde{P}_l^n - C_P \tilde{\vb*{v}}_r^n \vdot\tilde{ \vb*{v}}_r^n\frac{\rho_l}{\rho_g})}{\vb*{\varphi}_g}_{\Omega} \\
    & \quad + \inp*{\frac{1}{Re_g}\frac{\tilde{\grad}\alpha_g^n \vdot \tilde{\vb*{\tau}}_g^{n+\frac{1}{2}}}{\alpha_g^n}}{\vb*{\varphi}_g}_{\Omega} \\
    & \quad - \inp*{\frac{1}{Re_g}\tilde{\vb*{\tau}}_g^{n+\frac{1}{2}}}{\grad \vb*{\varphi}_g}_{\Omega} + \inp*{\frac{1}{Re_g} \vb*{n} \vdot \tilde{ \vb*{\tau}}_g^{n+\frac{1}{2}}}{\vb*{\varphi}_g}_{\Gamma_N}\\
    & \quad + \inp*{\frac{1}{Fr^2}  \tilde{\vb*{g}}}{\vb*{\varphi}_g}_{\Omega} -\inp*{\frac{3}{4}\frac{\rho_l}{\rho_g}\frac{C_{D}}{\tilde{d}_b}\norm{\tilde{\vb*{v}}_r^n}\tilde{\vb*{v}}_r^n}{\vb*{\varphi}_g}_{\Omega},
    \end{split}
    \end{align}%
\end{subequations}
where the subscript $\Omega$ denotes integral over the entire domain and $\Gamma_N$ denotes integral over boundaries where the Neumann boundary condition applies.
From \cref{eq:mom_NBC}, the normal component of the viscous stress tensor is equal to zero, substituting the relationship into the weak formulation:
\begin{subequations}
    \begin{align}
    \begin{split}
    \inp*{\frac{\tilde{\vb*{v}}_{l}^{*}- \tilde{\vb*{v}}_{l}^n}{\Delta t}}{\vb*{\varphi}_l}_{\Omega} + \inp*{\tilde{\vb*{v}}_{l}^n\vdot\tilde{\grad} \tilde{\vb*{v}}_{l}^n}{\vb*{\varphi}_l}_{\Omega} &= -\inp*{Eu_l \tilde{\grad} \tilde{P}_l^n}{\varphi}_{\Omega}  + \inp*{\frac{1}{Re_l}\frac{\tilde{\grad}\alpha_l^n \vdot \tilde{\vb*{\tau}}_l^{n+\frac{1}{2}}}{\alpha_l^n}}{\vb*{\varphi}_l}_{\Omega} \\
    & \quad - \inp*{\frac{1}{Re_l}\tilde{\vb*{\tau}}_l^{n+\frac{1}{2}}}{\grad \vb*{\varphi}_l}_{\Omega} + \inp*{\frac{1}{Fr^2}  \tilde{\vb*{g}}}{\vb*{\varphi}_l}_{\Omega}\\
    & \quad  +\inp*{\frac{3}{4}\frac{\alpha_g^n}{\alpha_l^n}\frac{C_{D}}{\tilde{d}_b}\norm{\tilde{\vb*{v}}_r^n}\tilde{\vb*{v}}_r^n}{\vb*{\varphi}_l}_{\Omega} \\
    & \quad -\inp*{ C_P \tilde{\vb*{v}}_r^n \vdot\tilde{ \vb*{v}}_r^n\frac{\tilde{\grad}\alpha_l^n}{\alpha_l^n}}{\vb*{\varphi}_l}_{\Omega},
    \end{split}\\
    \left.\tilde{\vb*{v}}_l^*\right|_{\Gamma_D} &= \tilde{\vb*{v}}_{l,BC}^{n+1},\\
    \begin{split}
    \inp*{\frac{\tilde{\vb*{v}}_{g}^{*}- \tilde{\vb*{v}}_{g}^n}{\Delta t}}{\vb*{\varphi}_g}_{\Omega} + \inp*{\tilde{\vb*{v}}_{g}^n\vdot\tilde{\grad} \tilde{\vb*{v}}_{g}^n}{\vb*{\varphi}_g}_{\Omega} &= -\inp*{Eu_g \tilde{\grad}\qty(\tilde{P}_l^n - C_P \tilde{\vb*{v}}_r^n \vdot\tilde{ \vb*{v}}_r^n\frac{\rho_l}{\rho_g})}{\vb*{\varphi}_g}_{\Omega} \\
    & \quad + \inp*{\frac{1}{Re_g}\frac{\tilde{\grad}\alpha_g^n \vdot \tilde{\vb*{\tau}}_g^{n+\frac{1}{2}}}{\alpha_g^n}}{\vb*{\varphi}_g}_{\Omega} \\
    & \quad - \inp*{\frac{1}{Re_g}\tilde{\vb*{\tau}}_g^{n+\frac{1}{2}}}{\grad \vb*{\varphi}_g}_{\Omega} + \inp*{\frac{1}{Fr^2}  \tilde{\vb*{g}}}{\vb*{\varphi}_g}_{\Omega}\\
    & \quad  -\inp*{\frac{3}{4}\frac{\rho_l}{\rho_g}\frac{C_{D}}{\tilde{d}_b}\norm{\tilde{\vb*{v}}_r^n}\tilde{\vb*{v}}_r^n}{\vb*{\varphi}_g}_{\Omega},
    \end{split}\\
    \left.\tilde{\vb*{v}}_g^*\right|_{\Gamma_D} &= \tilde{\vb*{v}}_{g,BC}^{n+1},
    \end{align}\label{eq:decoupled_IPCS_vten}%
\end{subequations}
where $\Gamma_D$ denotes the boundaries with Dirichlet boundary conditions.

\paragraph{Pressure}
The weak formulation of the pressure Poisson equation (\cref{eq:pp_discretized}) with the test function, $\varphi_p$ is:
\begin{equation}
\begin{split}
-\inp*{\sum_q Eu_q \alpha_q^n \tilde{ \grad} (\tilde{P}_l^{n+1}-\tilde{P}_l^n)}{\tilde{ \grad} \varphi_p}_{\Omega} &+ \inp*{\sum_q Eu_q \alpha_q^n \vb*{n} \vdot \tilde{ \grad} (\tilde{P}_l^{n+1}-\tilde{P}_l^n)}{\varphi_p}_{\Gamma_D} \\
&= \inp*{ \tilde{\div} \sum_q \qty(\frac{\alpha_q^n\tilde{\vb*{v}}_q^{*}}{\Delta t})}{\varphi_p}_{\Omega}.
\end{split}
\end{equation}
Given \cref{eq:pp_NBC}, the equation can be rewritten as:
\begin{subequations}
    \begin{align}
    -\inp*{\sum_q Eu_q \alpha_q^n \tilde{ \grad} (\tilde{P}_l^{n+1}-\tilde{P}_l^n)}{\tilde{ \grad} \varphi_p}_{\Omega} &= \inp*{ \tilde{\div} \sum_q \qty(\frac{\alpha_q^n\tilde{\vb*{v}}_q^{*}}{\Delta t})}{\varphi_p}_{\Omega},\\
    \left.\tilde{P}^{n+1}_l\right|_{\Gamma_N} &= \tilde{P}^{n+1}_{l,BC}.
    \end{align}\label{eq:decoupled_IPCS_p}%
\end{subequations}

\paragraph{Velocity Update}
The inner product of \cref{eq:vl_update_discretized,eq:vg_update_discretized} with $\vb*{\varphi}_q$ integrated over $\Omega$ is:
\begin{subequations}
    \begin{align}
    \inp*{\frac{\tilde{\vb*{v}}_{l}^{n+1}- \tilde{\vb*{v}}_{l}^*}{\Delta t}}{\vb*{\varphi}_l}_{\Omega} &=  -\inp*{ Eu_l \tilde{ \grad} (\tilde{P}_l^{n+1}-\tilde{P}_l^n)}{\vb*{\varphi}_l}_{\Omega},\\
    \left.\tilde{\vb*{v}}_l^{n+1}\right|_{\Gamma_D} &= \tilde{\vb*{v}}_{l,BC}^{n+1},\\
    \inp*{\frac{\tilde{\vb*{v}}_{g}^{n+1}- \tilde{\vb*{v}}_{g}^*}{\Delta t}}{\vb*{\varphi}_g}_{\Omega} &=  -\inp*{ Eu_g \tilde{ \grad} (\tilde{P}_l^{n+1}-\tilde{P}_l^n)}{\vb*{\varphi}_g}_{\Omega},\\
    \left.\tilde{\vb*{v}}_g^{n+1}\right|_{\Gamma_D} &= \tilde{\vb*{v}}_{g,BC}^{n+1}.
    \end{align}\label{eq:decoupled_IPCS_vcorr}%
\end{subequations}

\paragraph{Phase Fraction Update}
Lastly, the same procedure is repeated where the inner product of \cref{eq:alpha_update_discretized} with the test function $\varphi_\alpha$ is integrated over the simulation domain to give the following weak formulation:
\begin{subequations}
    \begin{align}
    \inp*{\frac{\alpha_g^{n+1} - \alpha_g^n}{\Delta t}}{\varphi_\alpha}_{\Omega} + \inp*{\tilde{\div}\qty(\alpha_g^{n+1}\tilde{\vb*{v}}_g^{n+1})}{\varphi_\alpha}_{\Omega} &= 0,\\
    \alpha_g^{n+1}|_{\Gamma_D} &= \alpha_{g,BC}^{n+1}.
    \end{align}\label{eq:decoupled_IPCS_alpha}%
\end{subequations}

\section{Grid Convergence}
\begin{figure}
    \centering
    \includegraphics[width=0.75\textwidth]{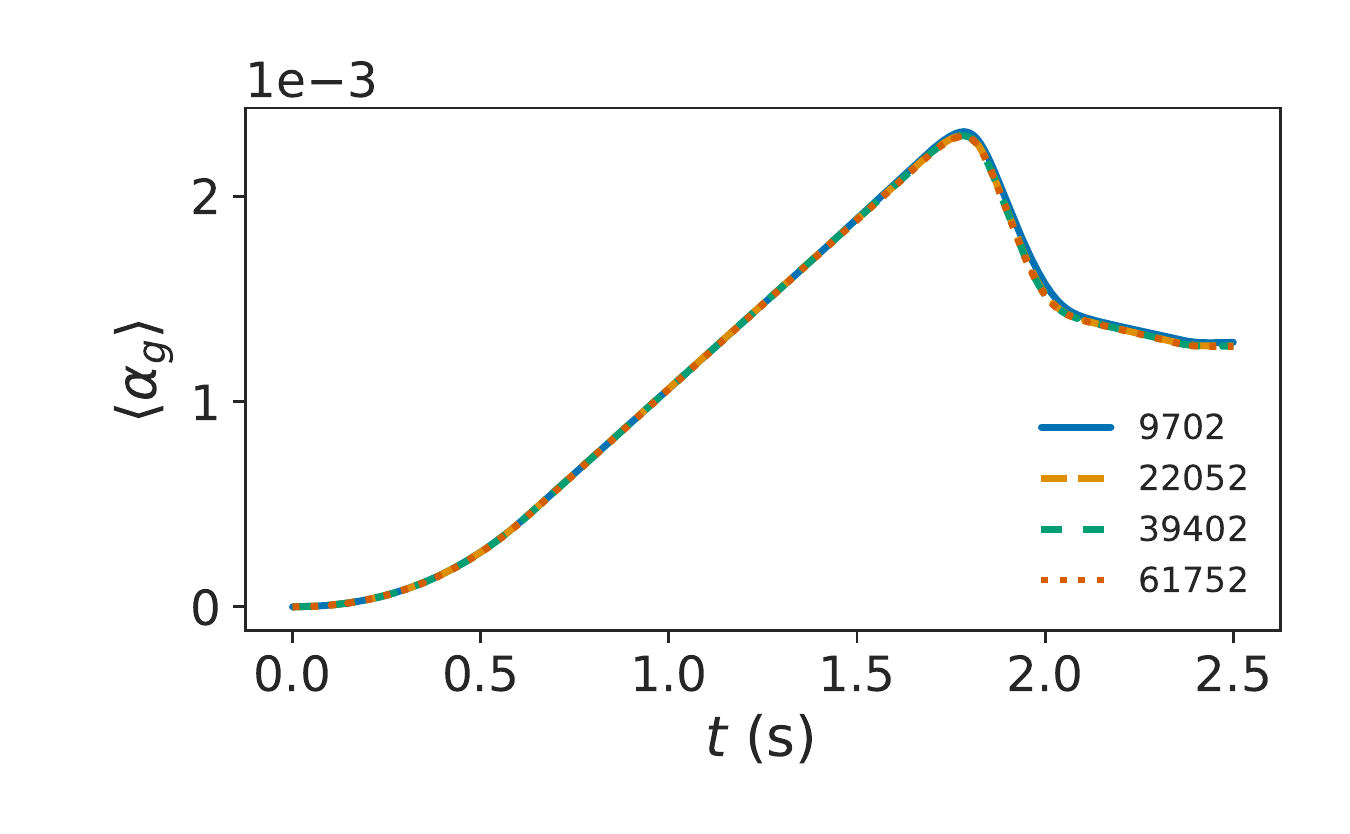}
    \caption{Time evolution of gas hold-up for different number of mesh elements.}
    \label{fig:results2_grid_convergence}
\end{figure}

To show mesh convergence, the bounded simulation with $P_c = P_d$ is repeated for three additional meshes, one coarser and two finer.
\Cref{fig:results2_grid_convergence} shows the time evolution of the gas hold-up inside in the simulation domain over the period of $\SI{2.5}{\second}$.
From \cref{fig:results2_grid_convergence}, the time evolution of the gas hold-up does not exhibit a discernible difference as the grid is refined from \num{22052} elements to \num{39402} elements.
Thus, the mesh with \num{22052} elements is used for all of the simulations reported in this study.

\begin{figure}
    \centering
    \includegraphics[width=0.75\textwidth]{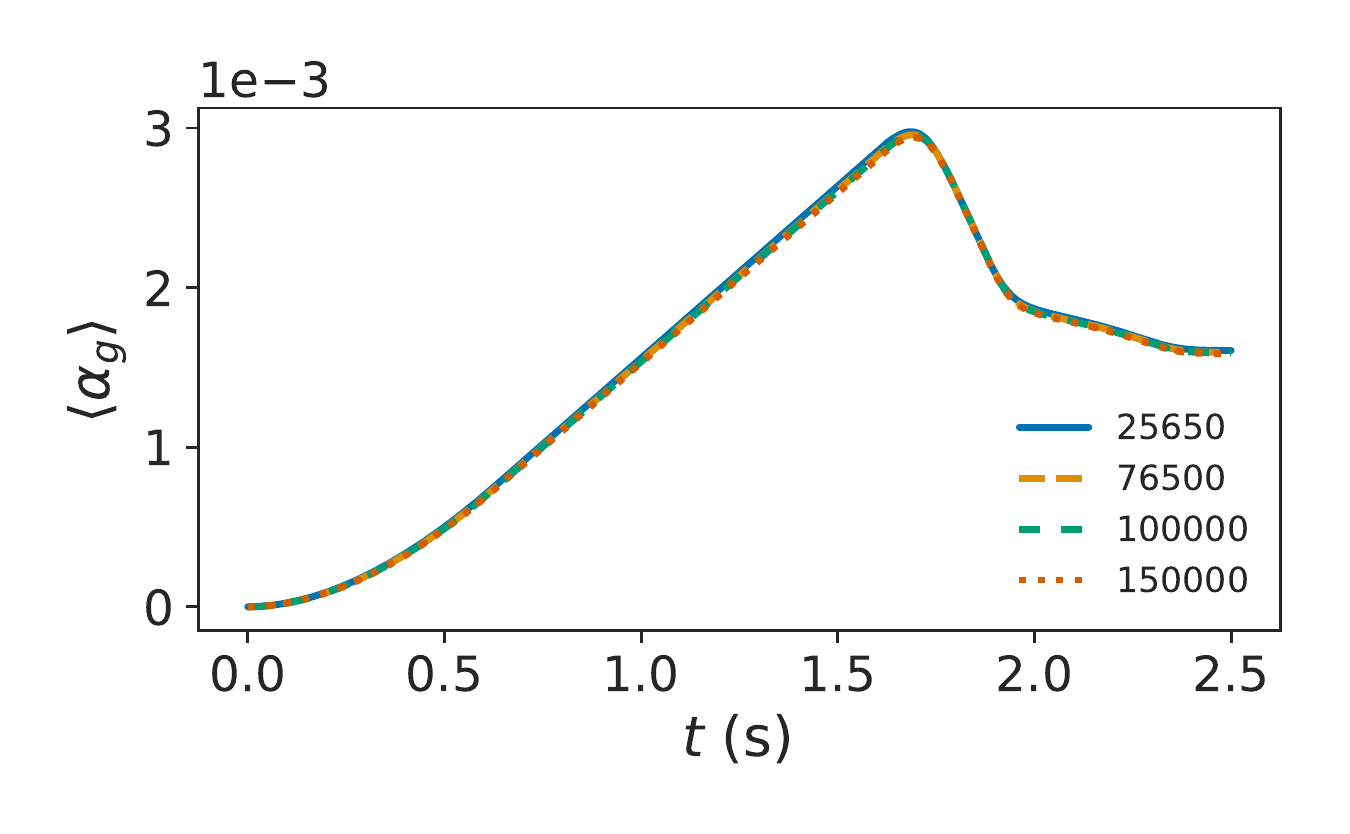}
    \caption{Time evolution of gas hold-up for different number of mesh elements from \texttt{OpenFOAM} simulations.}
    \label{fig:results2_grid_convergence_OF}
\end{figure}

Grid convergence in \texttt{OpenFOAM} is determined by performing simulations at increasing grid sizes, \numlist{25650;76500;100000;150000} mesh elements.
The time evolution of the gas hold-up is shown in \cref{fig:results2_grid_convergence_OF}.
The evolution of the gas hold-up does not significantly vary when the grid size is increased beyond \num{76500} and thus this grid size is used in this study.

\end{appendices}
\end{document}